\documentclass[11pt]{article}
\usepackage[total={18cm,21cm},top=2cm, left=2cm]{geometry}
\usepackage{amsmath, amssymb, amsfonts}
\usepackage{graphicx}
\usepackage{hyperref}
\usepackage{color}
\usepackage{epstopdf}
\usepackage{caption}
\usepackage{subcaption}
\newtheorem{theorem}{Theorem}
\newtheorem{proposition}{Proposition}
\newtheorem{lemma}{Lemma}
\newtheorem{corollary}{Corollary}

\begin{document}
\title{\Large\textbf{Numerical solution of the wave propagation problem in a plate.}}
\author{Manuel Cruz Rodriguez $^a$, Victoria Hern\'andez Mederos $^a$, Jorge Estrada Sarlabous $^a$,\\
Eduardo Moreno Hern\'andez $^a$, Ahmed Mansur Graver\'an $^a$ \\
$^a${\it Instituto de Cibern\'etica, Matem\'atica y F\'isica, ICIMAF, La Habana, Cuba}}

\date{}
\maketitle

\begin{abstract}
In this work, the propagation of an ultrasonic pulse in a thin plate is computed solving
the differential equations modeling this problem. To solve these equations finite differences
are used to discretize the temporal variable, while spacial variables are discretized using Finite
Element method. The variational formulation of the problem corresponding to a fixed value of time
is obtained and the existence an uniqueness of the solution is proved. The proposed approach leads
to a sequence of linear systems with the same sparse, symmetric and positive defined matrix.
The free software FreeFem++ is used to compute the approximated solution using polynomial triangular
elements. Numerical experiments show that velocities computed using the approximated displacements for
different frequency values are in good correspondence with analytical dispersion curves for the phase velocity.
\end{abstract}

\section{Introduction}
\label{section: Introduction}

Non-destructive tests (NDT) are used in the industry to inspect engineering structures that suffer
from the continuous influence of vibrations, loads and complicated environmental conditions, all producing
the appearance of defects, like corrosion, unbonds and delaminations. NDT play in consequence a very important role in modern
industry, where it is necessary to control manufacturing processes and reliability of products and systems,
whose premature failure would cause a hazard or economic loss. In the last decades, the development of NDT methods
has speed up due to several factors, such as the increase in manufacturing costs, the extreme service conditions of
components, and the design of structures with tighter adjustment of safety margins. Currently,  NDT techniques are applied
in energy and petrochemical industries, in the inspection of railway tracks, building structures, pipelines of gases and liquids
and in the study of aeronautical and aerospace structures, \cite{MAC00}, \cite{Galan02}, \cite{MAC03}.
%\cite{Perez}.

One of the main NDT techniques is based on application of ultrasonic guided waves (UGW) \cite{Hay03, Moreno y Col. 2015, Zhang18}.
Low frequency UGW propagate long distances in planar or tubular structures \cite{Low98}, because they employ the structure itself as a waveguide.
In consequence, UGW are suitable to inspect hidden or difficult to access areas, such as partially buried
structures covered with protective or insulating material or structures hidden behind other elements. In order to use
UGW in the industry it is necessary to study their propagation along the structures to be inspected and how they interact with possible
defects. This study is based on the computation of geometric {\it dispersion curves}, which describe the relationship between phase
(and therefore group) velocity and the frequency \cite{Hay03}. Dispersion curves enable the identification of the frequency intervals for which
waves propagate with less dispersion, that is the intervals where phase and group velocities have small variations. In this sense,
waves with almost constant phase velocity on a wide frequency interval are good candidates for NDT. Since dispersion curves can be
obtained analytically only for structures with simple geometry like bars and plates, simulation tools are frequently
used to compute them approximately.

Lamb waves are a special type of UGW propagating in solid plates. In the case of isotropic and homogeneous plates one can
find two types of Lamb waves: the symmetric and the anti-symmetric ones. Both modes are useful for Structural Health Monitoring (SHM) applications, where waves with comparable or greater with than the thickness of the plate and slight loss of amplitude magnitude are desirable, because they can be used to detect flaws at greater distances, see an excellent review on SHM in \cite{AbSh18}.

From the computational point of view, the solution of Lamb wave equations is very expensive, because  fine temporal and spacial
discretizations are necessary to reproduce the wave modes. The equations has been solved using several approaches, which depend on the
the boundary conditions, the length of the plate and other hypothesis. One of these approaches is the semi-analytical finite element method
(SAFEM) \cite{Galan03, Ahmad, Niel15, Her18}, based on the  assumption that the displacement of the wave may be factorized as the product of
an exponential function ( depending on time, wavenumber, angular frequency and the spatial variable defining the wave propagation axis)
and a function that only depends on spacial variables of the cross-section of the structure. SAFEM method leads to a family of
generalized eigenvalue problems depending on the wavenumber. The most common technique for solving Lamb wave equations is the
finite element method (FEM) \cite{Mace08, Sch11, Will12, Duc12} for spacial variables in combination with other methods for the
discretization of temporal variable. This approach leads to the solution of a sequence of linear systems. Recently, SAFEM and FEM have been
generalized to Isogeometric approach, which uses Non Uniform Rational B-splines as shape functions \cite{Her18, Will12, Duc12}.

%%%%%%%%%%%%%%%%%%%%%%%

\subsection{Related works}
Since Lamb waves are very important in NDT and Structural Health Monitoring (SHM), a lot of effort has been done to
study their propagation \cite{AbSh18, Gre83,Kos83,MAC00,Galan02,Galan03, MAC03, Ahmad,Mace08, Sch11, Will12, Duc12, Niel15,Her18}. Dispersion curves
of these waves have been calculated \cite{Achen73}, \cite{Rose99}
tracking the complex roots of the transcendent Rayleigh - Lamb frequency equation.
In \cite{Gomez-Ullate}, a study of Ultrasonic Systems based on multitrackers is carried out for the
detection of defects in plate-like structures. Moreover, a finite element model is proposed to study the propagation of Lamb
waves in airship components. Dispersion curves are also determined experimentally.
The numerical simulation of guided Lamb wave propagation in particle-reinforced composites is considered in
\cite {Weber y Col.}. The finite element method is used to perform parameter studies in order to better understand
how the propagation of the Lamb waves in these plates is affected by changes in the central frequency of the excitation
signal.
SAFEM is proposed in \cite{Galan02} to study Lamb wave scattering in homogeneous and sandwich isotropic
plates. An eigenvalue problem is solved to compute approximately the dispersion curves and the behavior of the error
of linear SAFEM for an infinite homogeneous layer of finite thickness is analized.
In \cite{Chen2012} a 2D FEM model is used to investigate the propagation of Lamb waves and its interaction
with a notch in a plate. FEM simulations can effectively evaluate the notch computing the delay time
between the reflection of the signal and first echo signal.
The Spectral Finite Element Method, the hierarchical p-FEM and the IgA approach are compared in \cite{Duc12} and \cite{Will12},
where they are used to compute the time-of-flight of Lamb waves propagating along a plate of finite length.
The study shows that high order finite element schemes are suitable to model the propagation of ultrasonic
Lamb waves. Moreover, the convergence properties of spectral finite elements and p-elements are very similar,
and IgA approach has the highest convergence rate.

\subsection{Our contribution}
In this paper we solve approximately the wave propagation equations on a isotropic plate under the effect of a pulse.
The corresponding system of two partial second order differential equations is solved using a combination of finite difference
and finite element methods. Finite difference approximation is used to discretize the time variable. For each fixed value of time,
Finite Element Method (FEM) is applied to the solution of the problem depending of spacial variables. The variational formulation
of the problem is obtained and the existence an uniqueness of the solution is proved. From the theoretical point of view the main contribution of this paper are the \textit{a priori} error estimates in two different norms for the approximated solution based on piecewise polynomials of degree $k$ (Theorems 2 and 3). This estimates show that the energy norm of error is proportional to $h^{k}$, while under an additional stability hypothesis the $L_2$ norm of the error is proportional to $h^{k+1}$.

Free software FreeFem++ \cite{Hecht15} is used to solve the spatial problem with polynomial triangular elements. Points on the phase velocity dispersion curve are computed approximately from the displacements obtained solving the differential problem. Numerical experiments confirm that the method produces good approximations of the dispersion curve and could be used for more complicated geometries.

\subsection{ Notation}

With $\boldsymbol{x}=(x,y)$, we use the following notation:
\begin{itemize}
\item[-] The classical scalar product of functions $u(\boldsymbol{x}),v(\boldsymbol{x}) \in L_2(\Omega)$  is defined as,
\begin{eqnarray}
\langle u(\boldsymbol{x}),v(\boldsymbol{x})\rangle_{L_2(\Omega)}&=&\iint\limits_{\Omega}u(\boldsymbol{x})v(\boldsymbol{x})\;d\boldsymbol{x}
\label{L2clasi}
%\|u(\boldsymbol{x})\| _{L_2(\Omega)}&=&\iint\limits_{\Omega}u(\boldsymbol{x})^2\;d\boldsymbol{x}
\end{eqnarray}
where $d\boldsymbol{x}=dxdy $.
\item[-] For order $n \times m$ matrix valued functions $\boldsymbol{A}(\boldsymbol{x})=(A_{ij}(\boldsymbol{x}))$ and $\boldsymbol{B}(\boldsymbol{x})=(b_{ij}(\boldsymbol{x}))$, with $a_{ij}(\boldsymbol{x}),b_{ij}(\boldsymbol{x}) \in L_2(\Omega)$, the scalar product is defined as
%we denote $\langle \cdot,\cdot \rangle_{L_2(\Omega)}$ the scalar product given by,
\begin{eqnarray}
\langle \boldsymbol{A}(\boldsymbol{x}),\boldsymbol{B}(\boldsymbol{x})\rangle_{L_2(\Omega)}&=& \sum _{i=1}^n\sum _{j=1}^m \, {\langle
a_{ij}(\boldsymbol{x}),b_{ij}(\boldsymbol{x}) \rangle}_{L_2(\Omega)} \label{peAB}
\end{eqnarray}
In particular, for $n=1$ and $m=2$ we obtain the following definition of scalar product between vector valued functions
 $\boldsymbol{u}(\boldsymbol{x})=(u_x(\boldsymbol{x}),u_y(\boldsymbol{x}))$ and
 $\boldsymbol{v}(\boldsymbol{x})=(v_x(\boldsymbol{x}),v_y(\boldsymbol{x}))$ in  $L_2(\Omega) \times L_2(\Omega)$,
\begin{eqnarray}
\langle \boldsymbol{u},\boldsymbol{v}\rangle_{L_2(\Omega)}= \iint\limits_{\Omega} \boldsymbol{u}(\boldsymbol{x})\cdot\boldsymbol{v}(\boldsymbol{x})\;d\boldsymbol{x}
&=&\iint\limits_{\Omega}(u_x(\boldsymbol{x})v_x(\boldsymbol{x})+u_y(\boldsymbol{x})v_y(\boldsymbol{x}))\;d\boldsymbol{x}
\label{peuv}
\end{eqnarray}
In both previous cases the corresponding norm is defined as usual in terms of the scalar product.
\item[-] The $\boldsymbol{:}$  product  between  matrices $\boldsymbol{A}(\boldsymbol{u})$  and  $ \boldsymbol{B}(\boldsymbol{v}) $
depending of vector value functions $\boldsymbol{u}$ and $\boldsymbol{v}$ is defined as \cite{Larson y Bengson},
\begin{equation}
\boldsymbol{A}\left( \boldsymbol{u}\right): \boldsymbol{B}\left( \boldsymbol{v} \right) = \sum \limits _{i=1}^{n}\sum \limits _{j=1}^{m}  a_{ij}\left( \boldsymbol{u} \right)b_{ij}\left(\boldsymbol{v}\right)
\label{A:B}
\end{equation}

\end{itemize}

%%%%%%%%%%%%%%%%%%%%%%%%%%%%%%%%%%%%%%%%%%%%%%%%%%%%%%%%
\section{The wave propagation problem on a plate}

\subsection{Formulation of the problem}
Lamb waves propagation in isotropic solid plates has been studied by many authors. It is well known that
the displacement of a particle is a vectorial function $\boldsymbol{u}(t,x,y)=(u_x(t,x,y),u_y(t,x,y))$ depending on
the temporal variable $t$, the spacial variables $\boldsymbol{x}=(x,y)$, the density of the material $\rho$ and the
Lam\'e constants $\lambda $ and $\mu $. The function $\boldsymbol{u}(t,x,y)$ is solution of the system of partial
differential equations \cite {Rose99},
\begin{equation}
\left \{ \begin{matrix}
\rho \dfrac{\partial ^2u_x}{\partial t^2} =\lambda\left(\dfrac{\partial ^{2}u_{x}}{\partial x^{2}} +\dfrac{\partial ^{2}u_{y}}{\partial x\partial y}\right)+2\mu \dfrac{\partial ^2 u_x}{\partial x^2}+\mu \left( \dfrac{\partial ^2 u_y}{\partial x \partial y}+\dfrac{\partial ^2u_x}{\partial y^2}\right) \\
\rho \dfrac{\partial ^2u_y}{\partial t^2}=\lambda \left(\dfrac{\partial ^2u_y}{\partial y^2}+\dfrac{\partial ^2u_x}{\partial x \partial y}\right) + 2\mu \dfrac{\partial ^2u_y}{\partial y^2} + \mu \left( \dfrac{\partial ^2u_x}{\partial x \partial y} +\dfrac{\partial ^2u_y}{\partial x^2} \right)
\end{matrix} \right.  \label{sistema}
\end{equation}

\begin{figure}[h!]
\centering
\includegraphics[trim={2cm 2cm 1cm 0.5cm},clip,scale=0.5]{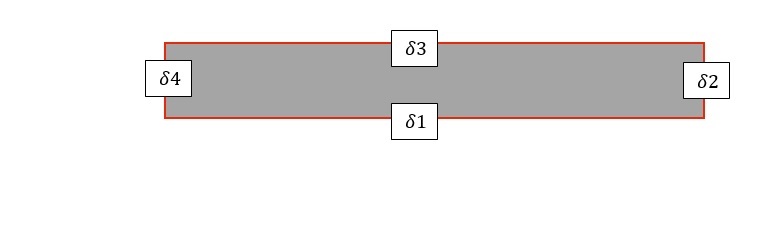}
\caption{Notation for the boundaries of the physical domain $ \Omega $.}\label{omega}
\end{figure}

Denote by $\boldsymbol{J}(\boldsymbol{u})$ the Jacobian matrix of the function $ \boldsymbol{u}(t,\boldsymbol{x})$ with respect to the spacial variables,
 $$\boldsymbol{J}(\boldsymbol{u}) = \begin{pmatrix}
 \dfrac{\partial u_x}{\partial x} & \dfrac{\partial u_x}{\partial y} \\
  \dfrac{\partial u_y}{\partial x} & \dfrac{\partial u_y}{\partial y}
 \end{pmatrix}
$$
In elasticity theory the strain tensor $\boldsymbol{S}(\boldsymbol{u})$ is defined as,
\begin{equation}
\boldsymbol{S}(\boldsymbol{u})=\begin{pmatrix}
S_{11}(\boldsymbol{u}) & S_{12}(\boldsymbol{u}) \\
S_{21}(\boldsymbol{u}) & S_{22}(\boldsymbol{u})
\end{pmatrix}=\dfrac{1}{2}(\boldsymbol{J}( \boldsymbol{u})+(\boldsymbol{J}( \boldsymbol{u}))^t)
= \begin{pmatrix}
\dfrac{\partial u_x}{\partial x} & \dfrac{1}{2} \left( \dfrac{\partial u_x}{\partial y} + \dfrac{\partial u_y}{\partial x} \right)\\
\dfrac{1}{2} \left( \dfrac{\partial u_x}{\partial y} + \dfrac{\partial u_y}{\partial x} \right) & \dfrac{\partial u_y}{\partial y}
\end{pmatrix} \label{S}
\end{equation}
Moreover, according to Hooke's law \cite{Rose99}, the stress tensor $\sigma(\boldsymbol{u})$ can be written
in terms of the strain as,
\begin{equation}
\sigma (\boldsymbol{u})=\begin{pmatrix}
\sigma_{11}(\boldsymbol{u}) & \sigma_{12}(\boldsymbol{u}) \\
\sigma_{21}(\boldsymbol{u}) & \sigma_{22}(\boldsymbol{u})
\end{pmatrix}=\begin{pmatrix}
\lambda tr(\boldsymbol{S}(\boldsymbol{u}))+2 \mu S_{11}{(\boldsymbol{u})} & 2\mu S_{12}{(\boldsymbol{u})} \\
2\mu S_{21}{(\boldsymbol{u})} & \lambda tr(\boldsymbol{S}(\boldsymbol{u}))+2 \mu S_{22}{(\boldsymbol{u})}
\end{pmatrix}
\label{T}
\end{equation}
where
\begin{equation}
tr(\boldsymbol{S}(\boldsymbol{u}))=\dfrac{\partial u_x}{\partial x}+\dfrac{\partial u_y}{\partial y} \label{tr}
\end{equation}
is the trace of $\boldsymbol{S}(\boldsymbol{u})$. Using the operator $\nabla=\left(\dfrac{\partial}{\partial x},\dfrac{\partial}{\partial y}\right)$, the divergence of the stress is defined as,
$$
\nabla\cdot \sigma =\begin{pmatrix}
\dfrac{\partial \sigma _{11}}{\partial x}+\dfrac{\partial \sigma _{21}}{\partial y} \\ \dfrac{\partial \sigma _{12}}{\partial x} + \dfrac{\partial \sigma _{22}}{\partial y}
\end{pmatrix}
$$
Computing partial derivatives of the components of the stress tensor, from the last expression and (\ref{T}) we obtain that $\nabla\cdot \sigma$ is the right hand side of equations (\ref{sistema}).
Hence, these equations can be written in compact form as,
\begin{equation}
\rho\frac{\partial ^{2}\boldsymbol{u}}{\partial t^2}=\nabla \cdot \sigma (\boldsymbol{u}) \label{rho}
\end{equation}
The physical domain $\Omega$ of our problem is the plate shown in the Figure \ref{omega}, with boundary $ \partial\Omega = \delta 1\cup\delta 2\cup\delta 3\cup\delta 4 $.

\begin{figure}[h!]
\centering
\includegraphics[scale=0.3]{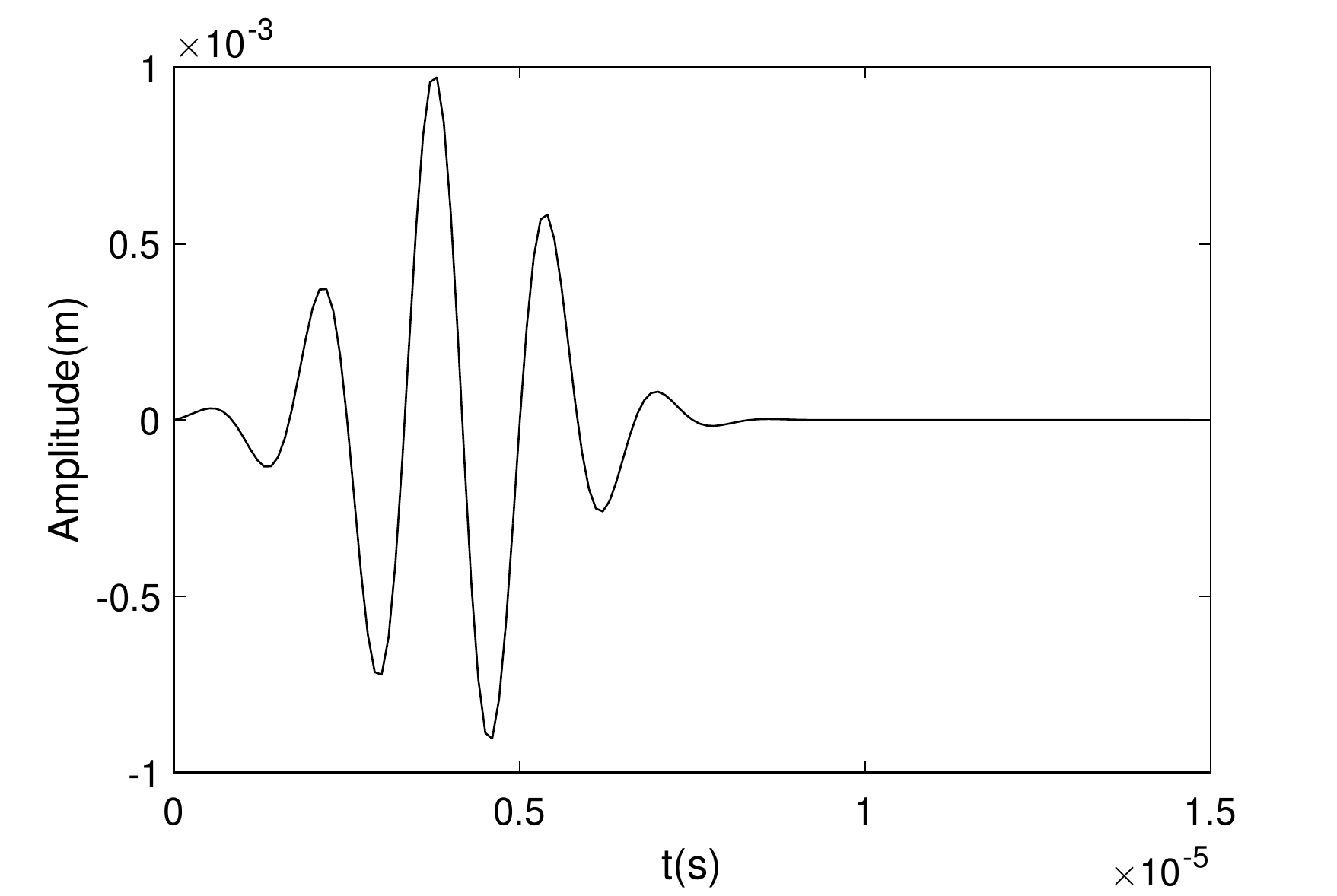}
\caption{Pulse $g(t)$ applied to the border $ \delta 4 $ of the plate. The parameters of the pulse are $f_0 = 600\ KHz$, $T_0=4.0 \cdot 10^{-6}\ s$,
$T=2.0 \cdot 10^{-6}\ s$, $\phi= 1.0 \cdot 10^{-3}\ m$ and $\alpha=1.1$.}
\label{pulso}
\end{figure}

In the boundary $ \delta4$, a pulse given by the function,
\begin{equation}
g(t)=\phi \sin (2\pi f_0 t)\exp{\left(-\alpha \frac{(t-T_0)^2}{T^2}\right)} \label{g}
\end{equation}
is applied for $t \in [0,t_f]$, with $t_f>0$ given. The parameters $f_0,\phi$ defining $g(t)$ are the frequency and the amplitude
of the pulse, while $T_0,\dfrac{T}{\sqrt{\alpha}}$ are the center and the width of the Gaussian factor of the pulse.
Figure \ref{pulso} shows a typical pulse $g(t)$.

In other words, in $\delta4$ the Dirichlet boundary condition,
\begin{equation}
\boldsymbol{u}(t,\boldsymbol{x})=(0,g(t)),\;\;\;\; \boldsymbol{x} \in \delta 4 \\
\label{Dirichlet}
\end{equation}
is imposed. Moreover, we assume that boundaries $ \delta 1, \delta 2 $ and $ \delta 3 $ are free, which means
that Neumann boundary conditions
\begin{equation}
\sigma (\boldsymbol{u}(t,\boldsymbol{x}))\cdot \boldsymbol{n}(\boldsymbol{x}) = \boldsymbol{0}, \;\;\;\;\boldsymbol{x} \in \delta 1, \delta 2, \delta 3
\label{Newman}
\end{equation}
hold, where $\boldsymbol {0} = (0,0) $ and $ \boldsymbol {n}(\boldsymbol{x}) $ is the vector normal to the boundary of the plate in the point
$\boldsymbol{x}$.
Finally, the initial conditions
\begin{equation}
\begin{matrix}
\boldsymbol{u}(0,\boldsymbol{x})=\boldsymbol{0}, & \boldsymbol{x}\in \Omega \\ \dfrac{\partial \boldsymbol{u}}{\partial t}(0,\boldsymbol{x})=\boldsymbol{0}, & \boldsymbol{x}\in \Omega
\end{matrix} \label{CI}
\end{equation}
are also imposed.

Summarizing, our mathematical problem is to calculate the function $\boldsymbol {u}(t,\boldsymbol {x})$ that satisfies the equations
(\ref {rho}) for $t\in [0,t_f]$ and $\boldsymbol{x}\in \Omega$, with the initial condition (\ref{CI}) and the boundary conditions (\ref{Dirichlet})and (\ref{Newman}). In the next sections we solve this problem using a combined method based, on finite differences
for time discretization and the classical finite element method for the discretization of spacial variables.

%%%%%%%%%%%%%%%%%%%%%%%%%%%%%%%%%%%%%%%%%%%%%%%%%%%%%%%%%%%%%%%%%%%%%%%%%%%%%%%%%%%%%%%%%%%%%%%%%%%%%%%%%%%%%%%%%%%%%

\subsection{Finite Differences discretization of time}

In order to discretize the temporal variable $t$ we build a {\it uniform mesh} in the interval $[0,t_f]$ with $N$ nodes,
$t_{i}=i\cdot\Delta t,\;\;\;\; i=0...N$, where $\Delta t=\frac{t_f}{N}$. In order to simplify the notation, from now on we use
\begin{equation}
\boldsymbol{u}_{i}(\boldsymbol{x})=\boldsymbol{u}(t_{i},\boldsymbol{x})
\label{variable}
\end{equation}
to denote the value of $\boldsymbol{u}(t,\boldsymbol{x})$ in the fixed time $t=t_i$.
The second derivative with respect to $t$ of the function $\boldsymbol{u}(t,\boldsymbol{x})$ in $t=t_i$ can be approximated using the
backward difference formula
\begin{equation}
\dfrac{\partial ^2 \boldsymbol{u}(t_{i},\boldsymbol{x})}{\partial t^2}=\dfrac{\boldsymbol{u}_{i}(\boldsymbol{x})- 2\boldsymbol{u}_{i-1}(\boldsymbol{x}) + \boldsymbol{u}_{i-2}(\boldsymbol{x})}{\left(\Delta t\right) ^2} \label{dudt}
\end{equation}
Substituting (\ref{dudt}) in (\ref{rho}) we obtain the discretization  with respect to time of the differential problem (\ref{rho}),
\begin{equation}
\begin{matrix}
\rho \boldsymbol{u}_{i}(\boldsymbol{x})-(\Delta t)^2 \left( \nabla \cdot \sigma (\boldsymbol{u}_{i}(\boldsymbol{x})) \right)=
\rho \left( 2 \boldsymbol{u}_{i-1}(\boldsymbol{x}) - \boldsymbol{u}_{i-2}(\boldsymbol{x}) \right), & \boldsymbol{x}\in \Omega ,
\end{matrix}\ \;\; i=1,...,N  \label{df2}
\end{equation}

In consequence, our problem is now reduced to find for each fixed time  $t=t_i,\;i=1,...,N$, the function $\boldsymbol{u}_i(\boldsymbol{x})$ that satisfies
the differential equations (\ref{df2}) with boundary conditions,
\begin{eqnarray}
\boldsymbol{u}_{i}(\boldsymbol{x})& = &\left(0,g\left(t_{i}\right)\right),\;\;\; \boldsymbol{x}\in \delta 4 \\
\sigma\left(\boldsymbol{u}_{i}(\boldsymbol{x})\right)\cdot \boldsymbol{n}(\boldsymbol{x}) & = & \boldsymbol{0}, \;\;\; \boldsymbol{x} \in \delta 1, \delta 2, \delta 3
\label{CF2}
\end{eqnarray}
From equations (\ref{df2}) it is clear that $\boldsymbol{u}_i(\boldsymbol{x})$ depends on the functions $ \boldsymbol{u}_{i-1}(\boldsymbol{x})$ and
$ \boldsymbol{u}_{i-2}(\boldsymbol{x})$ computed in the two previous time iterations. Taking into account (\ref{CI}), the initial conditions,
\begin{equation}
\boldsymbol{u}_{0}(\boldsymbol{x})=\boldsymbol{0},\;\;\;\boldsymbol{u}_{-1}(\boldsymbol{x})=\boldsymbol{0},\;\;\; \boldsymbol{x}\in \Omega
\label{CI3}
\end{equation}
are additionally imposed.

%%%%%%%%%%%%%%%%%%%%%%%%%%%%%%%%%%%%%%%%%%%%%%%%%%%%%%%%%%%%%%%%%%%%%%%%%%%%%%%%%%

\section{Finite Element discretization of spacial variables}

\subsection{Variational formulation}
To obtain the variational formulation, the problem (\ref {rho}) with boundary conditions \ref{CF2} is reduced to a problem with
homogeneous boundary conditions, writing  the solution as
\begin{equation}
\boldsymbol{u}_{i}(\boldsymbol{x})=\boldsymbol{u}^0_{i}(\boldsymbol{x})+\boldsymbol{u}_i^g(\boldsymbol{x}) \label{new u}
\end{equation}
where $ \boldsymbol {u}^0 _i(\boldsymbol{x}) $ is zero on $ \delta 4 $ and
\begin{equation}
\boldsymbol{u} _i^g(\boldsymbol{x})= \left\lbrace \begin{matrix} \boldsymbol{0}, & \boldsymbol{x} \in \Omega \setminus \delta 4 \\ \left(0,g(t_{i})\right), & \boldsymbol{x}\in \delta 4 \end{matrix} \right. \label{ug}
\end{equation}
Assuming that function $\boldsymbol{u}_i^g(\boldsymbol{x})$ is known, the problem is formulated now in terms of the
unknown is $ \boldsymbol {u}^0_ {i}(\boldsymbol{x})$. Replacing (\ref{new u}) in (\ref{df2}) we obtain the new problem to be solved
\begin{equation}
\rho \boldsymbol{u}^0_{i}(\boldsymbol{x})-(\Delta t)^2 \nabla \cdot \sigma\left(\boldsymbol{u}^0_{i}(\boldsymbol{x})\right)=\rho \left( 2 \boldsymbol{u}_{i-1}(\boldsymbol{x})
- \boldsymbol{u}_{i-2}(\boldsymbol{x}) \right)-\rho \boldsymbol{u}_i^g(\boldsymbol{x}) ,\;\;\; \boldsymbol{x}\in \Omega
\label{sistema2}
\end{equation}
with boundary conditions
\begin{eqnarray}
\boldsymbol{u}_{i}^0(\boldsymbol{x})&=&\boldsymbol{0},\;\;\; \boldsymbol{x}\in \delta 4 \\
\sigma\left(\boldsymbol{u}_{i}^0(\boldsymbol{x})\right)\cdot \boldsymbol{n}(\boldsymbol{x}) &=& \boldsymbol{0}, \;\;\; \boldsymbol{x} \in \delta 1, \delta 2, \delta 3
\label{CFs2}
\end{eqnarray}
and initial values
\begin{equation}
\boldsymbol{u}_{0}(\boldsymbol{x})=\boldsymbol{0},\;\;\;\boldsymbol{u}_{-1}(\boldsymbol{x})=\boldsymbol{0},\;\;\; \boldsymbol{x}\in \Omega
\label{CIs2}
\end{equation}

From now on,  and in order to simplify the notation, we suppress the dependence of functions on spacial variables $\boldsymbol{x}$.
Let $\mathcal{V}_0$ be the Hilbert space of functions
\begin{equation}
\mathcal{V}_0=\{\boldsymbol{v}\in H^2(\Omega)\times H^2(\Omega):\ \boldsymbol{v}|_{\delta 4}=\boldsymbol{0} \}
\label{V0}
\end{equation}
which consists of all vector valued functions $\boldsymbol{v} \in L_2(\Omega) \times L_2(\Omega)$ that possess weak and square-integrable
first and second derivatives and that vanish on $\delta 4$. Since functions in $\mathcal{V}_0$  vanish on a section of the boundary of $\Omega$, the scalar product $\langle \boldsymbol{u},\boldsymbol{v}\rangle_{\mathcal{V}_0}$ of functions $\boldsymbol{u}=(u_x,u_y)$ and $\boldsymbol{v}=(v_x,v_y)$ in $\mathcal{V}_0$
can be defined as \cite{Larson y Bengson},
\begin{equation}
\langle \boldsymbol{u},\boldsymbol{v}\rangle_{\mathcal{V}_0}=\iint\limits _{\Omega} (\nabla u_x \cdot \nabla v_x + \nabla u_y \cdot \nabla v_y) d\boldsymbol{x}
\label{dotV0}
\end{equation}
and the corresponding norm is
\begin{equation}
\left\| \boldsymbol{v} \right\| _{\mathcal{V}_{0}}^2 = \displaystyle \iint\limits _{\Omega} \left( \left( \dfrac{\partial v_x}{\partial x} \right) ^2 + \left( \dfrac{\partial v_y}{\partial x} \right) ^2 + \left( \dfrac{\partial v_x}{\partial y} \right) ^2 + \left( \dfrac{\partial v_y}{\partial y} \right) ^2 \right) d\boldsymbol{x}
\label{norma V0}
\end{equation}
Observe that according to (\ref{peAB}),
\begin{equation}
\left\| \boldsymbol{v} \right\| _{\mathcal{V}_{0}}=\left\| \boldsymbol{J}( \boldsymbol{v}) \right\| _{L_2(\Omega)}
\label{norma V0J}
\end{equation}

\medskip
The variational formulation of the problem is obtained computing the scalar product of both members of equation (\ref{sistema2})
with a function $ \boldsymbol{v} \in \mathcal{V}_0 $ and integrating the product over $ \Omega $,
\begin{equation}
\iint \limits _\Omega \rho \boldsymbol{u}^0_{i}(\boldsymbol{x})\cdot \boldsymbol{v}(\boldsymbol{x}) d\boldsymbol{x} -(\Delta t)^2 \iint \limits _\Omega \nabla \cdot \sigma (\boldsymbol{u}^0_{i}(\boldsymbol{x})\cdot\boldsymbol{v}(\boldsymbol{x}) d\boldsymbol{x}=\iint \limits _\Omega \rho \left( 2 \boldsymbol{u}_{i-1}(\boldsymbol{x}) - \boldsymbol{u}_{i-2}(\boldsymbol{x}) - \boldsymbol{u}_i^g(\boldsymbol{x}) \right)\boldsymbol{v}(\boldsymbol{x}) d\boldsymbol{x}
\label{var1}
\end{equation}
The second summand in the left hand side of (\ref{var1}) can be significantly simplified using Green's identity and taking into account boundary conditions (\ref{CFs2}) and the definition (\ref{ug}) of the function $ \boldsymbol{u}_i^g$. The resulting integral can be written in terms of the strain tensor as expressed in Lemma 1.

\begin{lemma}
\begin{equation}
-\iint \limits _\Omega \left(\nabla \cdot \sigma\right)\cdot\boldsymbol{v} d\boldsymbol{x}= 2\mu \iint \limits _\Omega \left(\boldsymbol{S}(\boldsymbol{u}_{i}^0):\boldsymbol{S}(\boldsymbol{v})\right)d\boldsymbol{x}+\lambda \iint\limits _\Omega (\nabla\cdot \boldsymbol{u}_{i}^0)(\nabla \cdot \boldsymbol{v})d\boldsymbol{x}
\label{2daint}
\end{equation}
where the $\boldsymbol{:}$  product is defined in \ref{A:B}.
\end{lemma}

Details of the proof of this Lemma can be found in \cite{Manuel}. Substituting (\ref{2daint}) in (\ref{var1}) we obtain the variational formulation of our problem given in the next proposition.

\begin{proposition} Variational formulation \\
The variational formulation of problem (\ref{sistema2}) is: to find $\boldsymbol{u}^0_{i} \in \mathcal{V}_0 $ such that for all $ \boldsymbol{v}\in \mathcal{V}_0$,
\begin{equation}
a(\boldsymbol{u}_{i}^0,\boldsymbol{v})=F(\boldsymbol{v})
\label{variacional}
\end{equation}
where
\begin{eqnarray}
a(\boldsymbol{u},\boldsymbol{v})&=&\rho \iint \limits _\Omega \left(\boldsymbol{u}\cdot \boldsymbol{v}\right) d\boldsymbol{x} +2\mu (\Delta t)^2 \iint \limits _\Omega
\left(\boldsymbol{S}(\boldsymbol{u}):\boldsymbol{S}(\boldsymbol{v})\right)d\boldsymbol{x}+\lambda (\Delta t)^2 \iint\limits _\Omega (\nabla\cdot \boldsymbol{u})(\nabla \cdot \boldsymbol{v})d\boldsymbol{x} \label{auv}\\
F(\boldsymbol{v})&=&\rho \iint \limits _\Omega \left( 2 \boldsymbol{u}_{i-1} - \boldsymbol{u}_{i-2} \right)\cdot\boldsymbol{v}\;d\boldsymbol{x}\label{Fv}
\end{eqnarray}
\end{proposition}

Observe that for $j=-1,0,...,i-1$ it holds $ u_j \in L_2(\Omega)\times L_2(\Omega)$, since $u_j=u_j^0+u_j^g$, with $ u_j^0 \in \mathcal{V}_0 \subset L_2(\Omega)\times L_2(\Omega)$ and $u_j^g \in L_2(\Omega)\times L_2(\Omega)$. In order to prove the existence and uniqueness of the solution of the variational problem (\ref{variacional}) in the next propositions we show some properties of $a(\boldsymbol{u},\boldsymbol{v})$ and $F(\boldsymbol{v})$.

\begin{proposition}\label{Propauv}
The functional $a(\boldsymbol{u},\boldsymbol{v}): \mathcal{V}_0 \times \mathcal{V}_0 \rightarrow \mathbb{R}$ given by (\ref{auv}),
\begin{itemize}
\item[a)] is a bilinear and continuous form
\item[b)] is coercive, i.e there is $C>0$ such that $a(\boldsymbol{u},\boldsymbol{u})\geq C \left\| \boldsymbol{u}\right\|_{\mathcal{V}_0}$
\end{itemize}
\end{proposition}

{\bf Proof}\\
a) The bilinearity of $a(\boldsymbol{u},\boldsymbol{v})$ is obtained immediately from the linearity of derivation and  integration and the linearity of the tensor $\boldsymbol{S}(\boldsymbol{u})$.
In order to prove the continuity observe that from (\ref{auv}), (\ref{A:B}), it holds
\begin{eqnarray*}
|a(\boldsymbol{u},\boldsymbol{v})| & \leq &  |\rho\iint\limits_{\Omega}\left(\boldsymbol{u}\cdot\boldsymbol{v}\right)d\boldsymbol{x}|+|2\mu (\Delta t)^2\iint\limits_{\Omega}\left(\boldsymbol{S}(\boldsymbol{u}):\boldsymbol{S}(\boldsymbol{v})\right)d\boldsymbol{x}| +|\lambda (\Delta t)^2\iint\limits_{\Omega}\left((\nabla \cdot\boldsymbol{u})(\nabla \cdot \boldsymbol{v})\right)d\boldsymbol{x}| \\
 & = &   \rho  |\langle \boldsymbol{u},\boldsymbol{v} \rangle _{L^2(\Omega)} | + 2\mu  (\Delta t)^2 | \langle \boldsymbol{S}(\boldsymbol{u}),\boldsymbol{S}(\boldsymbol{v})\rangle _{L^2(\Omega)} | +\lambda (\Delta t)^2   | \langle(\nabla \cdot\boldsymbol{u}),(\nabla \cdot \boldsymbol{v})\rangle _{L^2(\Omega)}|
\end{eqnarray*}
Hence, from the Cauchy-Schwartz inequality we obtain
\begin{equation}
|a(\boldsymbol{u},\boldsymbol{v})| \leq \rho\left\|\boldsymbol{u}\right\| _{L^2(\Omega)}\left\|\boldsymbol{v}\right\| _{L^2(\Omega)}+2\mu(\Delta t)^2\left\|\boldsymbol{S}(\boldsymbol{u})\right\|_{L^2(\Omega)}\left\| \boldsymbol{S}(\boldsymbol{v})\right\| _{L^2(\Omega)}
+\lambda(\Delta t)^2\left\|\nabla\cdot\boldsymbol{u}\right\|_{L^2(\Omega)}\left\|\nabla\cdot\boldsymbol{v}\right\|_{L^2(\Omega)}
\label{cont}
\end{equation}

The summands on the right hand side of (\ref{cont}) can be bounded in terms of $\left\| \boldsymbol{J}( \boldsymbol{u} ) \right\| _{L^2(\Omega )}$. In fact, from (\ref{S}) and (\ref{norma V0J}) we obtain that for all $\boldsymbol{u} \in \mathcal{V}_0 $,

\begin{equation}
\left\| \boldsymbol{S}(\boldsymbol{u})\right\| _{L^2(\Omega )}
%& = & \left\| \dfrac{1}{2}\left( \boldsymbol{J}(\boldsymbol{u})+(\boldsymbol{J}(\boldsymbol{u}))^t\right)\right\| _{L^2(\Omega )} \nonumber \\
  \leq  \dfrac{1}{2}\left( \left\| \boldsymbol{J}( \boldsymbol{u})  \right\| _{L^2( \Omega )} + \left\| \left(\boldsymbol{J}( \boldsymbol{u} ) \right)^t \right\| _{L^2( \Omega )} \right)  \\
  =  \left\| \boldsymbol{J}( \boldsymbol{u} ) \right\| _{L^2(\Omega )}=\left\| \boldsymbol{u} \right\| _{\mathcal{V}_0}
 \label{acot S}
\end{equation}
On the other hand,
\begin{eqnarray}
\left\|\nabla\cdot\boldsymbol{u}\right\|_{L^2(\Omega )}^2 & = &
%\left\| \dfrac{\partial u_x}{\partial x} +\dfrac{\partial u_y}{\partial y} \right\|_{L^2(\Omega )}=
\iint\limits_{\Omega}\left( \dfrac{\partial u_x}{\partial x} +\dfrac{\partial u_y}{\partial y} \right)^2 d\boldsymbol{x}
 \leq  \iint\limits_{\Omega}\left( \dfrac{\partial u_x}{\partial x} +\dfrac{\partial u_y}{\partial y} \right)^2 + \left( \dfrac{\partial u_x}{\partial x} -\dfrac{\partial u_y}{\partial y} \right)^2 d\boldsymbol{x} \nonumber \\
 &\leq &   2 \iint\limits_{\Omega}  \left( \dfrac{\partial u_x}{\partial x}\right)^2 + \left(\dfrac{\partial u_y}{\partial y} \right)^2 +\left(\dfrac{\partial u_x}{\partial y} \right)^2 + \left(\dfrac{\partial u_y}{\partial x} \right)^2 d\boldsymbol{x}
 =\left\| \boldsymbol{J}( \boldsymbol{u} )\right\|^2 _{L^2(\Omega )} \label{nablau}
\end{eqnarray}
Finally, from the inequality of Poincare-Friedrichs \cite{Martin}, it is easy to prove that for all $\boldsymbol{u} \in \mathcal{V}_0$ it holds,
\begin{equation}
 \|\boldsymbol{u}\|_{L^2(\Omega)} \leq \alpha \left\| \boldsymbol{J}( \boldsymbol{u}) \right\| _{L^2(\Omega )}
\label{Poincare}
\end{equation}
where $ \alpha>0 $ depends on $ \Omega $. Substituting (\ref{Poincare}),(\ref{nablau}) and (\ref{acot S}) in (\ref{cont}) we obtain,

\begin{equation}
|a(\boldsymbol{u},\boldsymbol{v})| \displaystyle\leq  C \left\| \boldsymbol{J}(\boldsymbol{u})\right\| _{L^2(\Omega)} \left\| \boldsymbol{J}(\boldsymbol{v})\right\| _{L^2(\Omega)} = C \left\| \boldsymbol{u}\right\|_{\mathcal{V}_0} \left\| \boldsymbol{v}\right\|_{\mathcal{V}_0}
\label{cotaabsauv}
\end{equation}
where $ C = \alpha^2  \rho  +2 \mu (\Delta t)^2 + 2 \lambda (\Delta t)^2 > 0 $. Hence, $a(\boldsymbol{u},\boldsymbol{v})$ is continuous \\

\medskip
b) to prove the coercivity of $a(\boldsymbol{u},\boldsymbol{v})$ observe that,
\begin{eqnarray*}
a(\boldsymbol{u},\boldsymbol{u}) & = & \displaystyle \rho\iint\limits_{\Omega}\left(\boldsymbol{u}\cdot\boldsymbol{u}\right) d\boldsymbol{x}+2\mu (\Delta t)^2\iint\limits_{\Omega}\left(\boldsymbol{S}(\boldsymbol{u}):\boldsymbol{S}(\boldsymbol{u})\right) d\boldsymbol{x} +\lambda (\Delta t)^2\iint\limits_{\Omega}\left(\nabla \cdot\boldsymbol{u}\right)^2 d\boldsymbol{x} \\
 & = & \displaystyle \rho\left\|\boldsymbol{u}\right\| _{L^2(\Omega)}^2 + 2\mu (\Delta t)^2 \left\| \boldsymbol{S}(\boldsymbol{u}) \right\| _{L^2(\Omega)}^2 + \lambda (\Delta t)^2 \left\| \nabla \cdot\boldsymbol{u}\right\| _{L^2(\Omega)}^2
\end{eqnarray*}
Moreover, since $ \rho, \mu, \lambda > 0 $ from the previous expression we obtain,
\begin{equation}
a(\boldsymbol{u},\boldsymbol{u})\geq 2\mu (\Delta t)^2  \left\| \boldsymbol{S}(\boldsymbol{u})\right\| _{L^2(\Omega)}^2
\label{coer}
\end{equation}
According to Korn inequality \cite{Brenner y Scott, Larson y Bengson} there is $\alpha >0$ depending on $\Omega$  and $\delta4$
such that for all $\boldsymbol{u} \in \mathcal{V}_0$ it holds,
\begin{equation}
\left\| \boldsymbol{S}(\boldsymbol{u})\right\| _{L^2(\Omega)}^2 \geq \alpha  \left\| \boldsymbol{J}(\boldsymbol{u})\right\| _{L^2(\Omega)}^2
\label{coer1}
\end{equation}
Hence, substituting (\ref{coer1}) in (\ref{coer}) and taking into account (\ref{norma V0J}) we obtain,
\begin{equation}
a(\boldsymbol{u},\boldsymbol{u})\geq C \left\| \boldsymbol{u}\right\|_{\mathcal{V}_0}^2
\label{coerci}
\end{equation}
where $ C=2 \alpha\mu (\Delta t)^2 >0$, i.e. $a(\boldsymbol{u},\boldsymbol{v})$ is coercive.
\medskip
$\Box$

\begin{proposition}\label{PropFv}
The functional $F(\boldsymbol{v}):\mathcal{V}_0 \rightarrow \mathbb{R}$ given by (\ref{Fv}) is  linear and continuous.
\end{proposition}

{\bf Proof}\\
The linearity of  $F(\boldsymbol{v})$ is obviously obtained from the linearity of the integral. \\
From (\ref{Fv}) and the Cauchy-Schwartz  inequality it holds for all $\boldsymbol{v}\in \mathcal{V}_0$
\begin{equation*}
\left| F(\boldsymbol{v}) \right| = \Big| \iint \limits _\Omega \left( \boldsymbol{f}_{i-1}\cdot\boldsymbol{v}\right)d\boldsymbol{x} \Big|
= \left| \langle \boldsymbol{f}_{i-1},\boldsymbol{v} \rangle _{L^2(\Omega)} \right|
\end{equation*}
where $\boldsymbol{f}_{i-1}=2 \boldsymbol{u}_{i-1} - \boldsymbol{u}_{i-2}$.
Moreover, using Cauchy-Schwarz inequality and substituting (\ref{Poincare}) in the last equation we obtain,

$$\left| F(\textbf{v}) \right|\leq  \alpha \| \boldsymbol{f}_{i} \| _{L^2(\Omega )} \left\| \boldsymbol{J}(\boldsymbol{v}) \right\| _{L^2(\Omega )} = C \|\boldsymbol{v} \| _{\mathcal{V}_0}$$
where $ C= \alpha \| \widehat{\textbf{f}}_{i}\| _{L^2(\Omega )}>0 $. Hence $F(\boldsymbol{v})$ is continuous.

\medskip
$\Box$

\begin{theorem}
The variational problem (\ref{variacional}) has a unique solution $\boldsymbol{u}^0_{i} \in \mathcal{V}_0 $.
\end{theorem}
{\bf Proof}\\
Propositions (\ref{Propauv}) and (\ref{PropFv}) guarantee that $a(\boldsymbol{u},\boldsymbol{v})$ and $F(\boldsymbol{v}) $ satisfy the hypothesis of Lax-Milgram Theorem \cite{Martin}. Hence the variational problem (\ref{variacional}) has a unique solution.
\medskip
$\Box$

%%%%%%%%%%%%%%%%%%%%%%%%%%%%%%%%%%%%%%%%%%%%%%%%%%%%%%%%%%%%%%%%%%%%%%%%%%%%%%%%%%%%%%%%%%%%%%%%%%%%%%%%%%%%%%%%%%

\subsection{The Galerkin approach}\label{section3.2}
The Galerkin method solves the variational problem (\ref{variacional}) looking for an approximated solution in a finite-dimensional subspace $\mathcal{V}$. In the classical FEM, this subspace is defined in terms of a partition of the physical domain $\Omega$ in a mesh of triangles or quadrilaterals. In this paper, we discretize $\Omega$ constructing a mesh of $N_t$ triangles $ \Delta _r $, $ r=1,...,N_t$ such that $\cup _{r=1}^{N_t} \Delta _r=\Omega$. The triangular mesh is denoted by $\tau^h$ and the finite-dimensional subspace of $\mathcal{V}$ is by $\mathcal{V}_h$, where $ h = \max \limits _{1 \leq r \leq N_t} h_r$ represents the {\it size of the triangulation}.\\

We choose to use a conforming finite element method, i.e., $\mathcal{V}_h$ is a subspace of the space $\mathcal{V}_0$, and the bilinear form and the linear form of the discrete problem are identical to the original ones, see \cite{Ciar1979}.  More precisely, the space $\mathcal{V}_h$ consists of global $C^0$ vector valued functions, that restricted to each triangle $ \Delta _r$ are polynomials of degree $k$ and vanish on $\delta 4$. More precisely, if we denote by $\mathbb{P}^k(\tau^h) $ the space of (scalar) piecewise polynomial functions on the triangulation $\tau^h$, then $\mathcal{V}_h$ is defined as

\begin{equation}
\mathcal{V}_h=\left\{ \boldsymbol{u}^h \in \mathbb{P}^k(\tau ^h) \times \mathbb{P}^k(\tau ^h) : \boldsymbol{u}^h|_{\delta 4}=\boldsymbol{0} \right\} \label{Vh}
\end{equation}
Denote by $\phi_r(\boldsymbol{x}),$ $r=1,...,N_n$ the basis of the space $ \mathbb{P}^k(\tau ^h)$. The number $N_n$ of basis functions depends on the degree $k$, for instance $N_n=N_v$ for $k=1$ and $N_n=N_v+N_e$ for $k=2$, where $N_v,N_e$ denote the number of vertices and edges of $\tau^h$ respectively. In general, $N_n$ is the number of {\it nodes} and if we denote by $\boldsymbol{x}_j,$ $j=1,...,N_n$ the $j$-th node, then the basis functions $\phi_r(\boldsymbol{x})$ are the classical {\it Lagrange functions} on $ \mathbb{P}^k(\tau ^h)$ defined by the property,
\begin{equation}
\phi_r(\boldsymbol{x}_j)=\left\{ \begin{matrix}
1,& r=j\\
0,& j\neq r
\end{matrix}\right. ,\;\;  j,r=1,...N_n
\label{phi}
\end{equation}
In consequence, if $ \boldsymbol{u} ^h \in \mathcal{V}_h $ it can be written as
\begin{equation}
\boldsymbol{u}^h(\boldsymbol{x})=\left(\sum\limits _{r=1}^{N_n} q_r^{(1)} \phi _r(\boldsymbol{x}),\sum\limits _{r=1}^{N_n} q_r^{(2)}\phi _r(\boldsymbol{x})\right)=\sum\limits _{r=1}^{N_n} q_r^{(1)} (\phi _r(\boldsymbol{x}),0)+\sum\limits _{r=1}^{N_n} q_r^{(2)} (0,\phi _r(\boldsymbol{x}))
\label{basexp}
\end{equation}
Observe that according to (\ref{phi}), $ \left(q_r^{(1)},q_r^{(2)} \right) = \boldsymbol{u}^h (\boldsymbol{x}_r),$ $r=1,...,N_n$. Moreover,
vector valued functions  $ (\phi_r, 0),(0, \phi_r), $ $r=1,...,N_n$ define a basis of $\mathcal{V}_h $. To simplify the expressions, we use the following notation for basis functions and coefficients,
\begin{eqnarray}
\left(\boldsymbol{\psi} _1,\boldsymbol{\psi} _2,...,\boldsymbol{\psi} _{2N_n}  \right) &=&\left( (\phi _1,0),(0,\phi _1),...,(\phi _{N_n},0),(0,\phi _{N_n}) \right) \label{Vh_base}\\
\left(d_1 , d_2 , ... ,d_{2N_n}  \right) &=&\left( q_1^{(1)} , q_1^{(2)} , ... , q_{N_n}^{(1)} , q_{N_n}^{(2)} \right)
\label{vectord}
\end{eqnarray}
With this notation, (\ref{basexp}) can be written as
\begin{equation}
\boldsymbol{u}^h(\boldsymbol{x})=\sum\limits _{r=1}^{2N_n} d_r  \boldsymbol{\psi} _r(\boldsymbol{x})
\label{uh_fin}
\end{equation}
The Galerkin approximation $\boldsymbol{u}^h$ is constructed demanding that it satisfies the variational formulation (\ref{variacional})
\begin{equation}
a(\boldsymbol{u}^h,\boldsymbol{v})=F(\boldsymbol{v}),\;\;for \;\; all \;\; \boldsymbol{v} \in \mathcal{V}_h
\label{formFEM}
\end{equation}
This is equivalent to require that,
\begin{equation}
a(\boldsymbol{u}^h,\boldsymbol{\psi}_j)=F(\boldsymbol{\psi}_j),\;\;j=1,...,2N_n
\label{Galerkin}
\end{equation}
because functions $\boldsymbol{\psi}_j,\;j=1,...,2N_n$ are a basis of $\mathcal{V}_h$.
Substituting (\ref{uh_fin}) in the last expression and taking into account the linearity of the integral and the strain $\boldsymbol{S}(\cdot)$ operators, we obtain a linear system of equations that can be written in matrix form as,

\medskip
\begin{equation}
\boldsymbol{A\,d}=\boldsymbol{b}
\label{sislin}
\end{equation}
where $\boldsymbol{d} =(d_1 , d_2 , ... ,d_{2N_n})^t$, $\boldsymbol{b}=(b_1,...,b_{2N_n})^t $ and $\boldsymbol{A}=(A_{rj}),\;r,j=1,...,2N_n$ with

\begin{eqnarray}
A_{rj}&=& \rho \iint \limits _\Omega \left(  \boldsymbol{\psi}_r \cdot  \boldsymbol{\psi} _j\right)d\boldsymbol{x} + 2\mu (\Delta t)^2 \iint \limits _\Omega \left( \boldsymbol{S}( \boldsymbol{\psi}_r):\boldsymbol{S}( \boldsymbol{\psi} _j)\right)d\boldsymbol{x} +\lambda (\Delta t)^2 \iint\limits _\Omega (\nabla\cdot  \boldsymbol{\psi}_r)(\nabla \cdot  \boldsymbol{\psi} _j)d\boldsymbol{x}
\label{Mat_Rig}\\
b_j&=&\iint \limits _\Omega \left( \boldsymbol{f}_{i-1} \cdot  \boldsymbol{\psi} _j\right) d\boldsymbol{x}
\end{eqnarray}

\begin{lemma}
Matrix $\boldsymbol{A}$ given by (\ref{Mat_Rig}) is sparse, symmetric and positive definite.
\end{lemma}

{\bf Proof}\\
Basis functions $\boldsymbol{\psi}_r,\;\;r=1,...,2N_n$ have compact support. Denote $supp(\boldsymbol{\psi}_r)=\{\boldsymbol{x} \in \Omega:\;\boldsymbol{\psi}_r(\boldsymbol{x})\neq 0 \}$. Then, from (\ref{S}) and (\ref{A:B}) it is clear that,
$$\boldsymbol{\psi}_r(\boldsymbol{x}) \cdot  \boldsymbol{\psi} _j(\boldsymbol{x})=\boldsymbol{S}( \boldsymbol{\psi}_r(\boldsymbol{x})):\boldsymbol{S}( \boldsymbol{\psi}_j(\boldsymbol{x}))=(\nabla\cdot  \boldsymbol{\psi}_r(\boldsymbol{x}))(\nabla \cdot  \boldsymbol{\psi} _j(\boldsymbol{x}))=0$$

for $\boldsymbol{x} \notin supp(\boldsymbol{\psi}_r) \bigcap supp(\boldsymbol{\psi}_j)$. It means that $A_{rj}=0$ if
$supp(\boldsymbol{\psi}_r) \bigcap supp(\boldsymbol{\psi}_j)=\emptyset$. \\
From the symmetry of the products $\cdot$ and $:$, the symmetry of $\boldsymbol{A}$ is obtained immediately.\\

Let $\boldsymbol{u}^h$ be a function in $\mathcal{V}_h$ different from $\boldsymbol{0}$. Since $\mathcal{V}_h$ is a subspace of $\mathcal{V}_0$ from (\ref{coerci}) we get,
\begin{equation}
0<\|\boldsymbol{u}^h\|_{\mathcal{V}_0} \leq C a(\boldsymbol{u}^h,\boldsymbol{u}^h)
\label{cota}
\end{equation}
with $C>0$. Substituting in (\ref{cota}) the expression (\ref{uh_fin}) of $\boldsymbol{u}^h$ in terms of the basis of $\mathcal{V}_h$ and using the bilinearity of $a(\cdot,\cdot)$ we obtain,
$$0< a\left(\sum\limits _{r=1}^{2N_n} d_r  \boldsymbol{\psi} _r(\boldsymbol{x}),\sum\limits _{s=1}^{2N_n} d_s  \boldsymbol{\psi} _s(\boldsymbol{x})\right)=\sum\limits _{r=1}^{2N_n} \sum\limits _{s=1}^{2N_n} d_r  d_s a(\boldsymbol{\psi} _r,\boldsymbol{\psi} _s) $$
The last expression can be written as $\boldsymbol{d}^t\boldsymbol{A}\boldsymbol{d}$. Since $\boldsymbol{u}^h\neq \boldsymbol{0}$ then $\boldsymbol{d}\neq \boldsymbol{0}$. Hence, we have proved that $\boldsymbol{d}^t\boldsymbol{A}\boldsymbol{d}>0$ for $\boldsymbol{d}\neq \boldsymbol{0}$, i.e $\boldsymbol{A}$ is positive definite.
\medskip
$\Box$

{\bf Remarks}
\begin{itemize}
\item As a consequence of the previous Lemma the linear system (\ref{sislin}) has a unique solution. Hence, the Galerkin approximation $\boldsymbol{u}^h$ exist and it is unique.
\item To compute the approximated solution $\boldsymbol{u}^h(\boldsymbol{x})$, for each fixed time $t=t_i,i=0,...,N$, we have to solve
the linear system (\ref{sislin}). The matrix $\boldsymbol{A}$ depends on the size $h$ of the triangulation and the time step $\Delta t$ and the right hand side $\boldsymbol{b}$ depends on time discretization.
\end{itemize}

\subsection{Convergence of the method}

In this section we study the {\it convergence} of the Galerkin approximation $\boldsymbol{u}^h$ when the size $h$ of the triangulation goes to zero to the solution $\boldsymbol{u}$ of the variational problem
\begin{equation}
a(\boldsymbol{u},\boldsymbol{v})=F(\boldsymbol{v})
\label{variationalf}
\end{equation}
where the bilinear form $a(\boldsymbol{u},\boldsymbol{v})$ is defined by (\ref{auv}) and
\begin{equation}
F(\boldsymbol{v})=\iint \limits _\Omega \boldsymbol{f}\cdot\boldsymbol{v}\;d\boldsymbol{x}\label{Fvf}
\end{equation}
for a given function $\boldsymbol{f}(\boldsymbol{u})=(f_x(\boldsymbol{u}),f_y(\boldsymbol{u})): \Omega \rightarrow \mathbb{R} \times \mathbb{R}$.\\

In order to measure the magnitude of the error $\boldsymbol{u}-\boldsymbol{u}^h$ we first consider in $\mathcal{V}_h$ the {\it energy scalar product} defined by the bilinear form (\ref{auv}),
\begin{equation}
\langle \boldsymbol{u} , \boldsymbol{v} \rangle_{E}=a(\boldsymbol{u},\boldsymbol{v})
\label{escproda}
\end{equation}
and the corresponding {\it energy norm}
\begin{equation}
\left\| \boldsymbol{u} \right\|_{E} =\sqrt{a(\boldsymbol{u},\boldsymbol{u})}
\label{normaa}
\end{equation}

\begin{lemma} Orthogonality \cite{{Brenner y Scott}}\\
% ver Brenner pag4
The error $\boldsymbol{u}-\boldsymbol{u}^h$ is {\it orthogonal} to $\mathcal{V}_h$ in the energy norm, i.e
\begin{equation}
\langle \boldsymbol{u}-\boldsymbol{u}^h, \boldsymbol{v} \rangle_{E}=a(\boldsymbol{u}-\boldsymbol{u}^h,\boldsymbol{v})=0
\label{ortog}
\end{equation}

for all $\boldsymbol{v}\in \mathcal{V}_h$.
\end{lemma}

From the previous Lemma, one can prove the next result.

\begin{lemma} Best approximation \cite{{Brenner y Scott}}\label{mejorap}\\
% ver Brenner pag5
Galerkin approximation $\boldsymbol{u}^h$ is the {\it best approximation} for $\boldsymbol{u}$  from $\mathcal{V}_h$ in the energy norm, i.e
\begin{equation}
\| \boldsymbol{u} - \boldsymbol{u} ^h \|_{E} = \min\limits _{\boldsymbol{v} \in \mathcal{V}_h} \| \boldsymbol{u} - \boldsymbol{v} \|_{E}
\label{aproxi}
\end{equation}
\end{lemma}

In the rest of this section we obtain bounds for the error $\boldsymbol{u} - \boldsymbol{u} ^h$, first in the  energy norm and later in the classic $L_2$ norm.  These results are obtained from standard estimates for interpolation error, see  for instance \cite{Ciar2013}.
\medskip

\begin{proposition} Polynomial interpolation error \label{acotacionpi}\\
Let $g:\Omega \subset \mathbb{R}^2 \rightarrow \mathbb{R}$ be a function in $H^r(\Omega)$ and $\tau^h$ a triangulation of $\Omega$ of size $h$. Denote by $\pi_kg$ the {\it piecewise polynomial} in  $\mathbb{P}^k(\tau^h)$, whose restriction to a triangle in $\tau^h$ with vertices $v_0,v_1$ and $v_2$, interpolates $g$ on the points
\begin{equation}
a_j:=\mu_0^jv_0+\mu_1^jv_1+\mu_2^jv_2,\;\;j=1,...,m_k
\label{intpoints}
\end{equation}
with $m_k=\frac{(k+1)(k+2)}{2}\;$, $\mu_0^j+\mu_1^j+\mu_2^j=1$ and
$\mu_0^j,\mu_1^j,\mu_2^j \in \left\{0,\frac{1}{k},...,\frac{k-1}{k},1 \right\}$.
Then, there is a constant $C>0$ (depending on the smallest angle in $\tau^h$ and on $k$) such that,
\begin{equation}
\| g - \pi_kg \|_1  \leq  C h^{q-1}\|g\|_{q}
\label{cotaerror}
\end{equation}
where $q=\min(k+1,r)$ and $\|g\|_{q}$ denotes the norm in the space $H^q(\Omega)$ given by
\begin{equation}
\|g\|_{q}=\left(\sum_{0 \leq i+j \leq q} \iint\limits _{\Omega} \left(\dfrac{\partial^i}{\partial x}\dfrac{\partial^j}{\partial y}\,g\right)^2\,d\boldsymbol{x}\right)^{1/2}
\label{normHq}
\end{equation}
\end{proposition}
The points (\ref{intpoints}) are known in the literature as {\it principal lattice of order $k$}.

\begin{corollary}
Under the same hypothesis of Proposition (\ref{acotacionpi}) it holds
\begin{equation}
\|D(g-\pi_kg)\|_{L^2(\Omega)} \leq C h^{q-1}\|g\|_{q}
\label{Derror}
\end{equation}
where $(Df)^2 = \left(\dfrac{\partial f}{\partial x} \right ) ^2 + \left ( \dfrac{\partial f}{\partial y} \right ) ^2$.
\end{corollary}

{\bf Proof}\\
Taking $q=1$ in (\ref{normHq}) we obtain that for any scalar function $f$
$$\|f\|_{1}^2=\iint\limits _{\Omega} f^2 + (Df)^2\,d\boldsymbol{x}=\|f\|_{L^2(\Omega)}^2+\|Df\|_{L^2(\Omega)}^2 $$
Hence,
\begin{equation}
\|Df\|_{L^2(\Omega)} \leq \|f\|_{1}
\label{DfL2f1}
\end{equation}
The bound (\ref{Derror}) is obtained applying the last inequality to $f:=g-\pi_kg$ and using (\ref{cotaerror}).
\medskip
$\Box$

\medskip

{\bf Remarks}\label{Rem}
\begin{itemize}
\item For a vector valued function $\boldsymbol{u}(\boldsymbol{x}) = (u_x(\boldsymbol{x}),u_y(\boldsymbol{x})$, the piecewise interpolation polynomial $\pi_k\boldsymbol{u}$ is the vector valued function $\pi_k\boldsymbol{u}=(\pi_ku_x,\pi_ku_y)$.
\item If $\boldsymbol{u}\in \mathcal{V}_0$ then  $\pi_1\boldsymbol{u} \in \mathcal{V}_h$, since $\pi_1\boldsymbol{u}$ is linear, interpolates $\boldsymbol{u}$ and $\boldsymbol{u}(\boldsymbol{x})=\boldsymbol{0}$ for $\boldsymbol{x} \in \delta4$.
\item For $k \geq 2$, $\pi_k\boldsymbol{u}$ is not necessarily null restricted to $\delta4$.
Nevertheless, in the rest of the section, and abusing of notation, we also denote by $\pi_k\boldsymbol{u}$ the piecewise polynomial in  $\mathbb{P}^k(\tau^h)\times \mathbb{P}^k(\tau^h)$ identically zero on $\delta4$ and interpolating $\boldsymbol{u}$ in the principal lattice points.
\item For the new function $\pi_k\boldsymbol{u}$, it  also holds the bound obtained applying (\ref{cotaerror}) to each component of the error function $\boldsymbol{e}:=\boldsymbol{u}-\pi_k\boldsymbol{u}$, since $\boldsymbol{e}$ vanishes on $\delta4$. Moreover, if $\boldsymbol{u}\in \mathcal{V}_0$ then $\pi_k\boldsymbol{u} \in \mathcal{V}_h$ for any $k \geq 1$.
\end{itemize}

\begin{proposition}Bound for the Jacobian of the interpolation error \label{acotacionJac}\\
Let $\boldsymbol{u}(\boldsymbol{x}) = (u_x(\boldsymbol{x}),u_y(\boldsymbol{x})):\Omega \rightarrow \mathbb{R} \times \mathbb{R}$ be a function in $H^q(\Omega) \times H^q(\Omega),\;q \geq 1$ and $\tau^h$ a triangulation of $\Omega$ of size $h$. Then, there are constants $C_1,C_2>0$ (depending on the smallest angle in $\tau^h$ and on $k$) such that,
\begin{equation}
\left\|\boldsymbol{J}(\boldsymbol{u} - \pi_k \boldsymbol{u}) \right\|_{L^2(\Omega)} \leq h^{q-1} \left(C_1 \left\|u_x\right\|_{q}^2 + C_2 \left\|u_y \right\|_{q}^2\right)^{1/2}
\label{cotaJac}
\end{equation}
\end{proposition}

{\bf Proof}\\
\begin{eqnarray}
\left\|\boldsymbol{J}(\boldsymbol{u} - \pi_1 \boldsymbol{u}) \right\| _{L^2(\Omega}^2
&=& \displaystyle \iint\limits _{\Omega}  \left ( \dfrac{\partial (u_x -\pi_1 u_x)}{\partial x} \right )^2 + \left ( \dfrac{\partial (u_x -\pi_1 u_x)}{\partial y} \right )^2 + \left ( \dfrac{\partial (u_y -\pi_1 u_y)}{\partial x} \right )^2 + \left ( \dfrac{\partial (u_y -\pi_1 u_y)}{\partial y} \right )^2 d\boldsymbol{x} \nonumber \\
%&=& \displaystyle \iint\limits _{\Delta _k} \left( \left | \dfrac{\partial (u_x -\pi u_x)}{\partial x} \right |^2 + \left | \dfrac{\partial (u_x -\pi u_x)}{\partial y} \right |^2 \right)dxdy + \displaystyle \iint\limits _{\Delta _k} \left( \left | \dfrac{\partial (u_y -\pi u_y)}{\partial x} \right |^2 + \left | \dfrac{\partial (u_y -\pi u_y)}{\partial y} \right |^2 \right)dxdy \nonumber \\
&=& \left\| D(u_x - \pi_k u_x) \right\| _{L^2(\Omega)}^2 + \left\|D(u_y - \pi_k u_y) \right\| _{L^2(\Omega)}^2 \label{des1}
\end{eqnarray}
Applying (\ref{DfL2f1}) with $f=u_x - \pi_k u_x$ in the first summand and $f=u_y - \pi_k u_y$ in the second summand we obtain,
$$\left\|\boldsymbol{J}(\boldsymbol{u} - \pi_1 \boldsymbol{u}) \right\| _{L^2(\Omega}^2 \leq
\left\| u_x - \pi_k u_x) \right\|_{1}^2+\left\| u_y - \pi_k u_y) \right\|_{1}^2$$
Using now the inequality (\ref{cotaerror}) it holds that for certain constants $C_1^{*},C_2^{*}>0$
$$\left\|\boldsymbol{J}(\boldsymbol{u} - \pi_k \boldsymbol{u}) \right\| _{L^2(\Omega)}^2 \leq
h^{2(q-1)}\left((C_1^{*})^2\| u_x \|_q^2+(C_2^{*})^2\| u_y\|_q^2\right)$$
The result (\ref{cotaJac}) with $C_1=(C_1^{*})^2$ and $C_2=(C_2^{*})^2$ is obtained from the last expression.
\medskip
$\Box$

The next theorem gives an \textit{a priori} estimation of the error $\boldsymbol{u}-\boldsymbol{u}^h$  in the {\it energy norm}.

\begin{theorem} A priori error estimate in energy norm \label{TeoerrorE}\\
Assume that the solution $\boldsymbol{u}:\Omega \rightarrow \mathbb{R} \times \mathbb{R}$ of (\ref{variationalf})
is a function in $H^{k+1}(\Omega) \times H^{k+1}(\Omega)$ . Given a triangulation of $\Omega$ of size $h$,  let $\boldsymbol{u}^h \in \mathcal{V}_h$ with $\mathcal{V}_h$ given by (\ref{Vh}) be the finite element solution defined by (\ref{formFEM}). Then,
\begin{equation}
\|\boldsymbol{u}- \boldsymbol{u}^h \|_{E} \leq h^k \left( \widetilde{C}_{1} \| u_x\|_{k+1}^2 + \widetilde{C}_{2} \| u_y\|_{k+1}^2\right)^{1/2}
\label{cotaerrorVh}
\end{equation}
where $\widetilde{C}_1,\widetilde{C}_2$ are positive constants.
\end{theorem}

{\bf Proof}\\
Recall that according to Lemma \ref{mejorap},  $\boldsymbol{u}^h$ is the best approximation for $\boldsymbol{u}$ from $\mathcal{V}_h$.
Since the piecewise polynomial $\pi_k\boldsymbol{u}$ interpolating $\boldsymbol{u}$ in the lattice points is in  $\mathcal{V}_h$  from  (\ref{aproxi}) we obtain,
\begin{equation}
\| \boldsymbol{u} - \boldsymbol{u} ^h \|_{E} \leq \| \boldsymbol{u} - \pi_k \boldsymbol{u} \|_{E}
\label{acotacionPi}
\end{equation}
Moreover, from (\ref{normaa}) and (\ref{cotaabsauv}) we get,
\begin{eqnarray}
\| \boldsymbol{u} - \pi_k \boldsymbol{u}\|_{E}^2&=&a(\boldsymbol{u} - \pi_k \boldsymbol{u},\boldsymbol{u} - \pi_k \boldsymbol{u}) \nonumber\\
&\leq& C\left \|\boldsymbol{J}(\boldsymbol{u}- \pi_k\boldsymbol{u})\right \| _{L^2(\Omega)}^2
\label{cotaJ}
\end{eqnarray}
with $C=\rho \alpha ^2 + 2\mu (\Delta t)^2 + 2\lambda (\Delta t)^2 >0 $.
Finally,  since $\boldsymbol{u} \in H^{k+1}(\Omega)  \times H^{k+1}(\Omega)$ we can use the bound (\ref{cotaJac}) obtained in Proposition \ref{acotacionJac} for the Jacobian of the interpolation error with $q=k+1$. Taking square root in (\ref{cotaJ}) and substituting this bound it follows,
$$\| \boldsymbol{u} - \pi_k \boldsymbol{u}\|_{E} \leq \sqrt{C}\,h^k \left(C_1 \|u_x\|_{k+1}^2 + C_2 \|  u_y\| _{k+1}^2\right)^{1/2} $$
The result (\ref{cotaerrorVh}) is obtained substituting the last expression in (\ref{acotacionPi}) with $\widetilde{C}_1=CC_1$ and $\widetilde{C}_2=CC_2$.
\medskip
$\Box$

The energy norm is useful since it allows to obtain a simple derivation of the \textit{a priori} error estimate (\ref{cotaerrorVh}). However, it is not a natural norm such as the $L_2$-norm. In the next theorem we show that under an additional {\it stability condition} it can be shown that the $L_2$ error of the FEM solution is of order $h^{k+1}$, if we use polynomials of degree $k$ to construct the approximated solution.
In the proof we follow an approach based on the duality principle and the stability condition similar to other authors  \cite{Stewart98},
\cite{Larson y Bengson}, \cite{Brenner y Scott}.

\begin{theorem}\label{Theo3} A priori error estimate in $L_2$ norm\\
Assume that the solution $\boldsymbol{u}=(u_x(\boldsymbol{x}),u_y(\boldsymbol{x})):\Omega \rightarrow \mathbb{R} \times \mathbb{R}$ of (\ref{variationalf}) is a function in $H^{k+1}(\Omega) \times H^{k+1}(\Omega)$ that satisfies the stability condition
\begin{equation}
\left(\|u_x\|_{2}^2+\|u_y\|_{2}^2\right)^{1/2} \leq \widetilde{C}\left(\|f_x\|_{L_2(\Omega)}^2+\|f_y\|_{L_2(\Omega)}^2\right)^{1/2}
\label{stahipo}
\end{equation}
where $\widetilde{C}$ is a positive constant. Given a triangulation of $\Omega$ of size $h$, let $\boldsymbol{u}^h \in \mathcal{V}_h$ with $\mathcal{V}_h$ given by (\ref{Vh}) be the finite element solution defined by (\ref{formFEM}). Then, the following \textit{a priori} error estimate in the $L_2$ norm holds,
\begin{equation}
\|\boldsymbol{u}- \boldsymbol{u}^h \| _{L_2(\Omega)} \leq h^{k+1} \left( \widehat{C}_{1} \| u_x\| _{k+1}^2 + \widehat{C}_{2}
\|u_y\| _{k+1}^2\right)^{1/2}
\end{equation}
where $\widehat{C}_1,\widehat{C}_2$ are positive constants.
\end{theorem}

{\bf Proof}\\
The proof makes use of the duality principle. Let $\boldsymbol{w}=(w_x,w_y)$ be the solution of the dual problem
\begin{eqnarray}
\rho\boldsymbol{w}-(\Delta t)^2\nabla\cdot\sigma(\boldsymbol{w})&=&\boldsymbol{e} \label{ecuacionerror}\\
\boldsymbol{w}(\boldsymbol{x})&=&\boldsymbol{0},\;\;\boldsymbol{x} \in \delta4 \label{boundcond}
\end{eqnarray}
where $\boldsymbol{e}=\boldsymbol{u}- \boldsymbol{u}^h$. After  scalar multiplication of both members of equation (\ref{ecuacionerror}) by $\boldsymbol{e}$ and integrating we obtain,
\begin{equation}
\rho\iint\limits _{\Omega} \boldsymbol{w}\cdot \boldsymbol{e}\;d\boldsymbol{x} - (\Delta t)^2\iint\limits _{\Omega}\nabla\cdot\sigma(\boldsymbol{w})\cdot \boldsymbol{e}\;d\boldsymbol{x}= \iint\limits _{\Omega} \boldsymbol{e}^2 \;d\boldsymbol{x}
\end{equation}
Substituting (\ref{2daint}) in the second summand it holds,
\begin{equation}
\rho\iint\limits _{\Omega} \boldsymbol{w}\cdot \boldsymbol{e}\;d\boldsymbol{x} +2\mu (\Delta t)^2 \iint\limits _{\Omega} \boldsymbol{S}(\boldsymbol{w}):\boldsymbol{S}(\boldsymbol{e})\;d\boldsymbol{x} + \lambda(\Delta t)^2 \iint\limits _{\Omega} (\nabla \cdot \boldsymbol{w})(\nabla \cdot \boldsymbol{e})\;d\boldsymbol{x} = \iint\limits _{\Omega} \boldsymbol{e}^2 \;d\boldsymbol{x}
\end{equation}

Using the scalar product (\ref{escproda}) and the corresponding energy norm, from the previous equation we obtain
\begin{equation}
\|\boldsymbol{e} \|_{L^2(\Omega)}^2=a(\boldsymbol{e},\boldsymbol{w})
\label{eL2a}
\end{equation}
Since the piecewise interpolating polynomial $\pi_k\boldsymbol{w}$ satisfies
$\pi_k\boldsymbol{w}(\boldsymbol{x})=\boldsymbol{0}$ for $\boldsymbol{x} \in \delta4$ (see Remark \ref{Rem}) it is a function in
$\mathcal{V}_h$. Hence, according to (\ref{ortog}) $a(\boldsymbol{e},\pi_k\boldsymbol{w})=0$ and we can write
(\ref{eL2a}) as,
\begin{eqnarray}
\|\boldsymbol{e}\|_{L^2(\Omega)}^2&=&a(\boldsymbol{e},\boldsymbol{w})-a(\boldsymbol{e},\pi_k\boldsymbol{w})\nonumber\\
                              &=&a(\boldsymbol{e},\boldsymbol{w}-\pi_k\boldsymbol{w})=\langle \boldsymbol{e},\boldsymbol{w}-\pi_k\boldsymbol{w} \rangle_{E} \nonumber\\
                              &=&\|\boldsymbol{e}\|_{E} \;\|\boldsymbol{w}-\pi_k\boldsymbol{w}\|_{E}  \label{cotaeL2}
\end{eqnarray}
where we have used the definition (\ref{escproda}) of the scalar product and the Cauchy-Schwartz inequality.\\
\medskip
Now observe that the solution $\boldsymbol{w}$ of the dual problem is a function in $H^2(\Omega) \times H^2(\Omega)$ ( since the right hand side of the dual problem is the error $\boldsymbol{e}$ which lies in $H^1(\Omega) \times H^1(\Omega)$). Hence, from (\ref{cotaJ}) and the bound (\ref{cotaJac}) obtained in Proposition \ref{acotacionJac} for the Jacobian of the interpolation error with $q=2$
we obtain
\begin{eqnarray}
\|\boldsymbol{w}-\pi_k\boldsymbol{w}\|_{E} &\leq&
                                   C \left \|\boldsymbol{J}(\boldsymbol{w}- \pi_k\boldsymbol{w})\right\|_{L^2(\Omega)}\nonumber\\
                            &\leq& C h \left(C_1 \| w_x\|_{2}^2 + C_2 \| w_y\|_{2}^2\right)^{1/2}\nonumber\\
                            &\leq& \widehat{C} h \left(\| w_x\|_{2}^2 + \| w_y\|_{2}^2\right)^{1/2}
\label{cotaerroePikw}
\end{eqnarray}
where $\widehat{C}=C\max\{C_1,C_2\}$. Furthermore, applying the stability hypothesis (\ref{stahipo}) to the dual problem solution $\boldsymbol{w}$ with right hand side $\boldsymbol{e}=(e_x,e_y)$ we obtain
$$(\|w_x\|_{2}^2+\|w_y\|_{2}^2)^{1/2} \leq \widetilde{C} (\|e_x\|_{L_2(\Omega)}^2+\|e_y\|_{L_2(\Omega)}^2)^{1/2}=\widetilde{C}\|\boldsymbol{e}\|_{L_2(\Omega)}$$
for certain constant $\widetilde{C}>0$. Hence, substituting the last inequality in (\ref{cotaerroePikw}) it holds,
\begin{equation}
\|\boldsymbol{w}-\pi_k\boldsymbol{w}\|_{E} \leq \widehat{C}\widetilde{C}h\|\boldsymbol{e}\|_{L_2(\Omega)}
\label{cotaerroePikwnew}
\end{equation}

Finally, substituting in (\ref{cotaeL2}) the inequality (\ref{cotaerroePikwnew}) and the error estimate (\ref{cotaerrorVh}) in energy norm given in Theorem \ref{TeoerrorE},  we obtain
\begin{equation}
\|\boldsymbol{e}\|_{L^2(\Omega)}^2 \leq h^{k+1} \widehat{C}\widetilde{C}\left( \widetilde{C}_{1} \| u_x\|_{k+1}^2 + \widetilde{C}_{2}
\| u_y\|_{k+1}^2\right)^{1/2}\|\boldsymbol{e}\|_{L_2(\Omega)}
\label{cotacasifinal}
\end{equation}

The result is obtained from this inequality with $\widehat{C}_1=(\widehat{C}\widetilde{C})^2\widetilde{C}_{1}$ and
$\widehat{C}_2=(\widehat{C}\widetilde{C})^2\widetilde{C}_{2}$.
\medskip
$\Box$

\section{Numerical results and discussion}

To solve numerically the wave propagation problem we wrote a FreeFem++ \cite{Hecht15} code that computes the approximated solution
using either \textit{linear} or \textit{quadratic} Lagrange finite elements. The code first constructs a triangulation $\tau$ of the plate and then, for each time step, it solves the corresponding linear system (\ref{sislin}) using the direct sparse solver multi-frontal method UMFPACK.

In wave propagation modelling it is customary to impose bounds for  mesh size $h$ and time step $\Delta t$. More precisely, $h$ is selected as a fraction of the wavelength $\Lambda(f_0)=\frac{c_L}{f_0}$, where $c_L=C(f_0)$ is the phase velocity of the wave, and
$ \Delta t$ is chosen such that  $\Delta t \leq {\Delta t}_{cr}=\frac{h}{c_L}$ (\textit{CFL condition}), where ${\Delta t}_{cr}$
is the transit time of the wave through the smallest element in the model \cite{CFL28},\cite{Bathe96}, \cite{Drozdz}. Consequently, in the following numerical experiments we are using $h < \frac{\Lambda}{4}$ and  $\Delta t <\frac{1}{4f_0}$, such that both upper bounds are  satisfied.

\begin{figure}[b]
\centering
\includegraphics[scale=0.4]{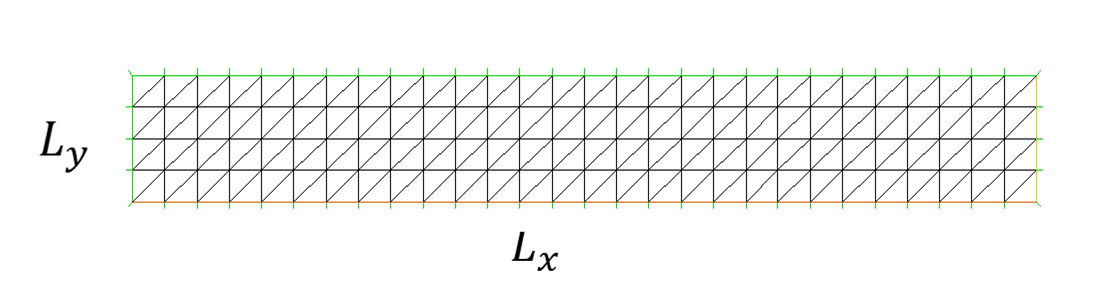}
\caption{Quasi-uniform and shape regular triangulation  $\tau$ fitted to $\Omega$.}
\label{Fig:triangulacion}
\end{figure}

Given a set of vertices on the boundary $\partial \Omega$ of the plate,  FreeFem++  constructs the triangulation $\tau$  containing these boundary vertices using the Delaunay-Voronoi algorithm. In all our experiments, $ny$  uniformly distributed vertices are generated on each boundaries $\delta2$ and $\delta4$. The corresponding amount $nx$ of uniformly distributed  vertices  on each boundary $\delta1$ and $\delta3$ is computed as $nx=\lceil ny(L_x/L_y)\rceil$. In this way, for any value of $ny$, the triangles of the corresponding mesh $\tau$ are as close to isosceles as possible ( see Figure \ref{Fig:triangulacion}) and consequently $\tau$ is a quasi-uniform and shape regular triangulation fitted to $\Omega$ \cite{Ciar1979}.

In our experiments we work with an aluminium plate with density $\rho=2700\,Kg/m^3$, Poisson ratio $P=0.334 $, Young's modulus $E=7.0\cdot 10^{10}\,N/m^2$ and  Lam\'e constants
$\mu=2.624\cdot 10^{10}\,N/m^2$ and $\lambda=5.279\cdot 10^{10}\,N/m^2$.  The size of the plate in the directions $x$ (wave propagation) and $y$ are $L_x=5.0\cdot 10^{-2}\,m$ and $L_y=1.0\cdot 10^{-3}\,m$ respectively.

All numerical results of this section are obtained in a PC with i5 processor and $6 $ Gb of RAM. Moreover, spacial variables are measured in meters and time is measured in seconds.

\subsection{Phase velocity dispersion curve}\label{seccionPhasevel}
In dispersive media the phase velocity of the wave depends on the frequency. This dependence is described by
the phase velocity dispersion curves, which are very important in the industry, since they  enable the identification of the frequency intervals for which waves propagate with less dispersion. In the case of thin plates, the  phase velocity dispersion curve $d(f_0)$ is the parametric curve,
$$d(f_0)=\left(\frac{L_y}{\Lambda(f_0)},\frac{C(f_0)}{c_0} \right)$$
where $c_0=\sqrt{\frac{E}{\rho}}$ is a constant, $L_y$ is the width of the plate and  $C(f_0)$ and  $\Lambda(f_0)=\frac{C(f_0)}{f_0}$ are  the phase velocity  and the wavelength corresponding to the frequency $f_0$.

\begin{table}[h!]
\centering
\caption{Parameters of the pulse applied at the boundary $ \delta 4 $.}
\begin{tabular}{|c|c|c|c|c|c|c|c|} \hline
$\boldsymbol{f_0}$ & $\boldsymbol{100\, KHz \,}$& $\boldsymbol{200\, KHz \,}$& $\boldsymbol{300\, KHz \,}$& $\boldsymbol{600\, KHz \,}$& $\boldsymbol{700\, KHz \, }$& $\boldsymbol{900\, KHz \,}$& $\boldsymbol{1100\, KHz \,}$ \\ \hline
$\alpha$ & $2.5\cdot 10^{-2}$& $1.5\cdot 10^{-1}$& $8\cdot 10^{-2}$& $1.1$& $1.5$ & $2$& $3.2$\\ \hline
$\phi$ & $10^{-3}\ m$& $10^{-3}\ m$& $10^{-3}\ m$& $10^{-3}\ m$& $10^{-3}\ m$ &$10^{-3}\ m $&$10^{-3}\ m $\\ \hline
$T_0$ & $23 \cdot 10^{-6}\ \,s$& $12 \cdot 10^{-6}\ \,s$& $8 \cdot 10^{-6}\ \,s$& $4 \cdot 10^{-6}\ \,s$& $3 \cdot 10^{-6}\ \,s$& $2.6 \cdot 10^{-6}\ \,s$&$2.3 \cdot 10^{-6}\ \,s$ \\ \hline
$T$ & $2 \cdot 10^{-6}\ \,s$& $2 \cdot 10^{-6}\ \,s$& $1 \cdot 10^{-6}\ \,s$& $2 \cdot 10^{-6}\ \,s$& $2 \cdot 10^{-6}\ \,s$& $2 \cdot 10^{-6}\ \,s$& $2 \cdot 10^{-6}\ \,s$ \\ \hline
\end{tabular}
\label{tablaf}
\end{table}

\begin{figure}[h!]
\centering
\includegraphics[scale=0.4]{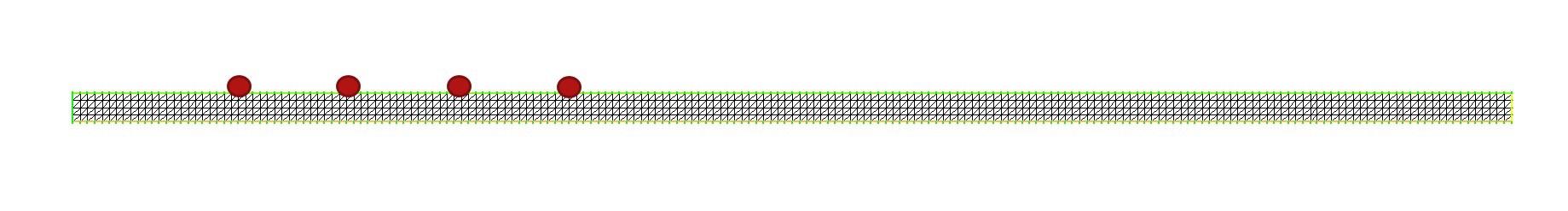}
\caption{Selected points $\boldsymbol{p}_i=(x_i,L_y),\;i=1,...,4$ on the top of the plate with $x_1=1.0\cdot 10^{-2},\,$ $x_2=1.3\cdot 10^{-2},\,$ $x_3=1.6\cdot 10^{-2}$ and $x_4=1.9\cdot 10^{-2}$.}
\label{Fig:puntos}
\end{figure}

\begin{figure}[h!]
\centering
\includegraphics[scale=0.4]{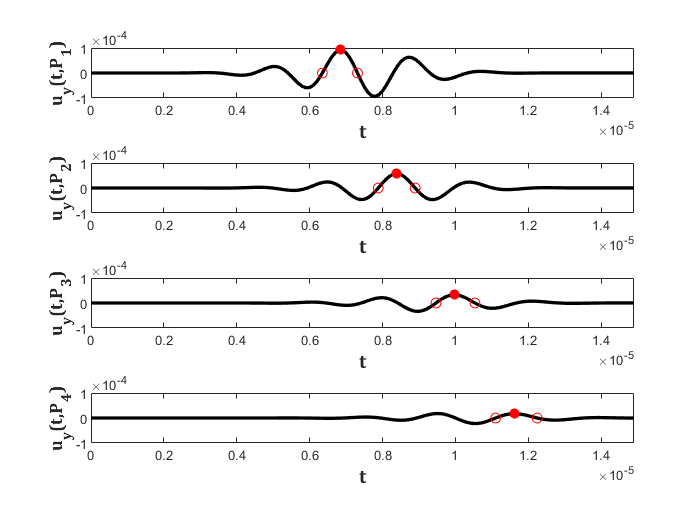}\hspace{2cm}
\includegraphics[scale=0.35]{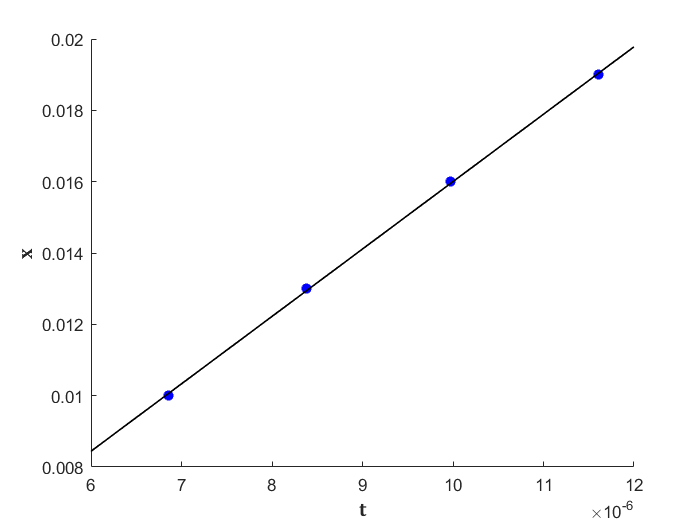}
\caption{Left: Graphic of the displacement in the vertical direction  $u_y(t,\boldsymbol{p}_i)$ of the selected points $\boldsymbol{p}_i=(x_i,L_y),\;i=1,...,4$. Red bullets: point on $u_y(t,\boldsymbol{p}_i)$ corresponding to the arrival time $\widetilde{t}_i$ of the pulse at the point $\boldsymbol{p}_i, \;i=1,...,4$, for $f_0= 600\,KHz$. Red circles: two zeros of the function $u_y(t,\mathbf{p}_i)$. Right: Plot of  arrival times $\widetilde{t}_i$ versus $x_i$, for $i=1,...,4$.}
\label{Fig:curvastemporales}
\end{figure}

In this section we compute some  points on the phase velocity dispersion curve $d(f_0)$ proceeding as follows. First, we select
4 fixed points  on the top of the plate with coordinates $\mathbf{p}_i=(x_i,L_y),\;i=1,...,4$ with
$x_1=1.0\cdot 10^{-2},\,$ $x_2=1.3\cdot 10^{-2},\,$ $x_3=1.6\cdot 10^{-2}$ and $x_4=1.9\cdot 10^{-2}$, see Figure \ref{Fig:puntos}.
For each of the pulses with parameters in  Table \ref{tablaf} (depending on $f_0$), we compute the numerical solution of the wave propagation problem using  \textit{quadratic} FEM with a mesh of size $h= 3.53 \cdot 10^{-4}$. The problem is solved for $t \in [0,1.5 \cdot 10^{-5}]$, using finite differences for a sequence $t_j,\;j=0,1...,150$ of equidistant values of time, with time step $\Delta t = 1.0\cdot 10^{-7}$. In Figure \ref{Fig:curvastemporales} left we show, for each fixed point $\mathbf{p}_i,\;i=1,...,4$, the curve interpolating the sequence of vertical displacements $u_y^h(t_j,\mathbf{p}_i)$ obtained for the values of $t_j,\;j=0,...,150$. The arrival time $\widetilde{t}_i$ of the pulse at the point $\mathbf{p}_i,\;i=1,...,4$ is computed as the value of $t$ for the second maximum
( represented in Figure \ref{Fig:curvastemporales} as a bullet) of the function $u_y^h(t,\mathbf{p}_i)$. Moreover, the zeros $t_i^{min},t_i^{max}$ (circles in Figure \ref{Fig:curvastemporales}) of the function $u_y(t,\mathbf{p}_i)$ closest to $\widetilde{t}_i$ are also computed, with $t_i^{min}< \widetilde{t}_i <t_i^{max}$ for $i=1,...,4$. These values are used in next two sections to compute a suitable interval to compare the approximated solution with Lamb waves and also to compare the FEM solutions obtained for meshes of decreasing size. It holds that for a fixed  value of $f_0$ the points $(\widetilde{t}_i,x_i),\;i=1,...,4$ are approximately on a line; the slope of this line is the phase velocity $C(f_0)$ corresponding to the frequency $f_0$. Figure \ref{Fig:curvastemporales} right shows the points $(\widetilde{t}_i,x_i),\;i=1,...,4$ and the fitting line for $f_0= 600\,KHz$.

\begin{figure}[hb]
\centering
\includegraphics[scale=0.4]{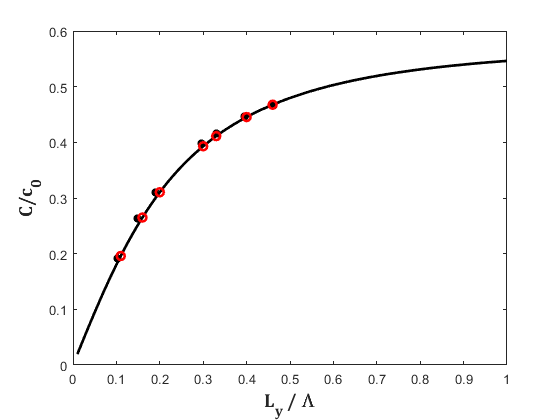}
\caption{Theoretical  phase velocity dispersion curve (continuous line), points ($\bullet$)  computed with Lagrange  \textit{quadratic} FEM for all frequency values $f_0$ from Table \ref{tablaf} and orthogonal projections ( red circles) on the dispersion curve.}
\label{Fig:dispersion}
\end{figure}

Repeating the previous methodology for the frequencies $f_0$ in Table \ref{tablaf}, we compute the phase velocity $C(f_0)$ corresponding to each frequency $f_0$. In Figure \ref{Fig:dispersion}  we show  the points $\left(\frac{L_y}{\Lambda(f_0)},\frac{C(f_0)}{c_0} \right)$  computed from  the numerical solution. For comparison, we also show the graph of the theoretical phase velocity dispersion curve (continuous line), defined by an implicit curve $F\left(\frac{L_y}{\Lambda(f_0)},\frac{C(f_0)}{c_0} \right)=0$, see \cite{Rose99}, \cite{Martincek},\cite{Moreno y Col. 2015}. As we can see, the points computed from the FEM solution are very close to the theoretical dispersion curve; in fact the distance from these points to their orthogonal projections is of order $ 1.0\cdot 10^{-2}$. Hence, the previous strategy could be useful to compute approximately the phase velocity dispersion curve for more complicated geometries or nonelastic and anisotropic materials, where the theoretical phase velocity dispersion curves are unknown \cite{MAC00}, \cite{MAC03}.

\subsection{Qualitative study of FEM approximation, graphic comparison}\label{secqual}

\begin{figure}[h!]
\centering
\includegraphics[scale=0.18]{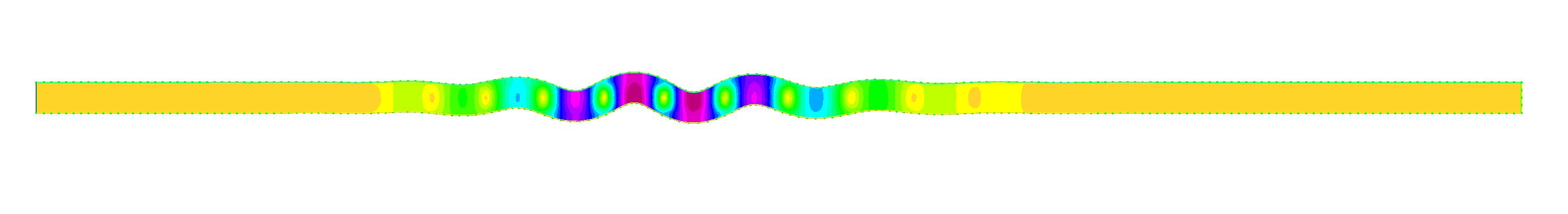} \vspace{-0.5cm}\\
\includegraphics[scale=0.18]{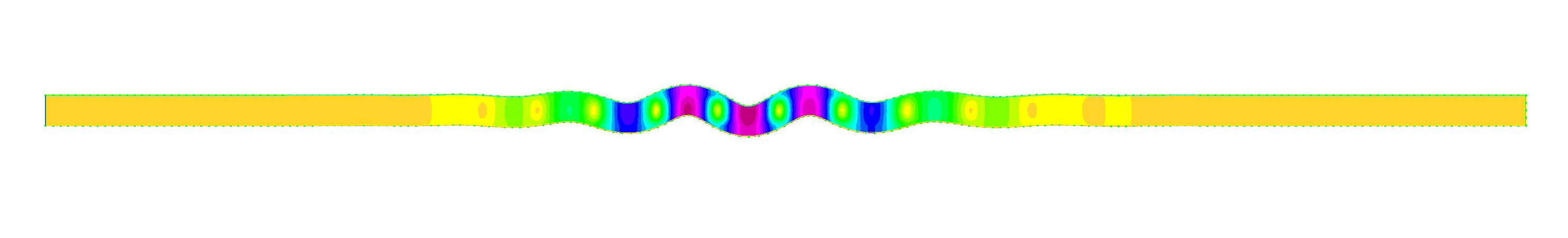} \vspace{-0.5cm}\\
\includegraphics [scale=0.18]{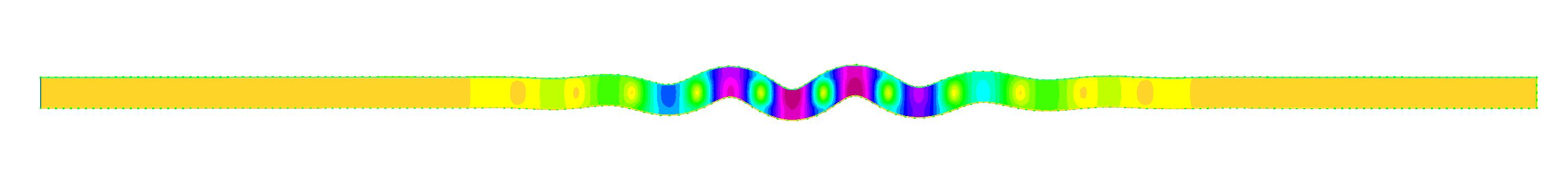}
\caption{Deformation of the plate after emitting a pulse on the boundary $\delta 4$. From top to bottom graphics for $t=1.23 \cdot 10^{-5}$,$t=1.32 \cdot 10^{-5}$ and $t=1.40 \cdot 10^{-5}$. Colors correspond to the intensity of the deformation (the norm of the displacement vector), yellow means no displacement, pink represents the maximum displacement.}
\label{Fig:desplaza}
\end{figure}

For  the pulse given by (\ref{g}) with parameter $f_0=600 \ KHz$ in Table \ref{tablaf}, we compute in this section the approximated solution $\mathbf{u}^h$ using finite differences with  $\Delta t = 1.0 \cdot 10^{-7}$ and Lagrange  \textit{quadratic} finite elements, with a mesh of size $ h = 1.41 \cdot 10^{-4}$. Figure \ref{Fig:desplaza} shows the deformation of the plate in three different times after emitting the pulse on the boundary $\delta 4$. Colors in this figure correspond to the $L_2$ norm of the displacement vector.\\

The Dirichlet boundary condition (\ref{Dirichlet}) associated with the pulses in (\ref{g}) produces antisymmetric displacements. Certainly our wave propagation problem does not have an exact solution, but we may check if our approximated displacement $\boldsymbol{u}^h$ has a sound behavior (from the  point of view of the physics) comparing it to the exact solution of a similar  wave propagation problem. Hence, we compare $\mathbf{u}^h$ with the antisymmetric solution $\mathbf{u}$ of the Lamb wave equations (\ref{rho}), for an \textit{infinite} plate in the direction $x$ of the wave propagation (and length $L_y$ in the $y$ direction), assuming Neumann boundary conditions. The exact solution  $\boldsymbol{u}$ of this problem is given by $\mathbf{u}(t,x,y)=Re(\mathbf{g}(y)e^{i(\omega t -kx)})$,
where  $\omega=2\pi f_0$ is the angular frequency, $k=2\pi/\Lambda$ is the wavenumber, $\Lambda$ is the wavelength and $\mathbf{g}(y)=(g_x(y),g_y(y))$ is given in \cite{Martincek}, \cite{Rose99}.

\begin{figure}[h!]
\centering
\includegraphics[scale=0.25]{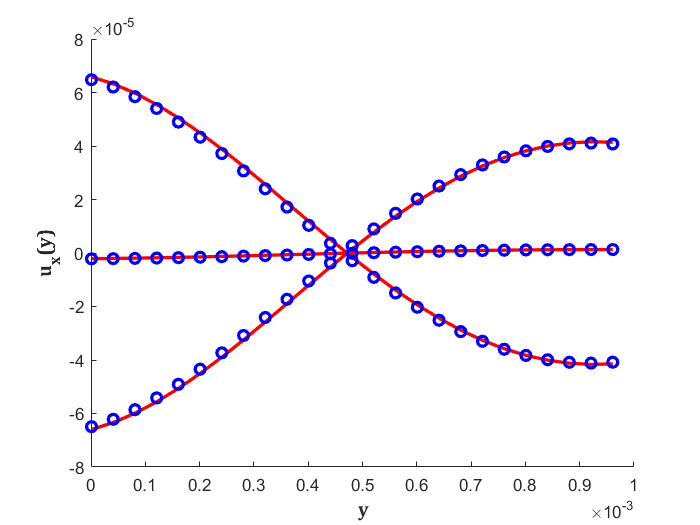}
\includegraphics[scale=0.25]{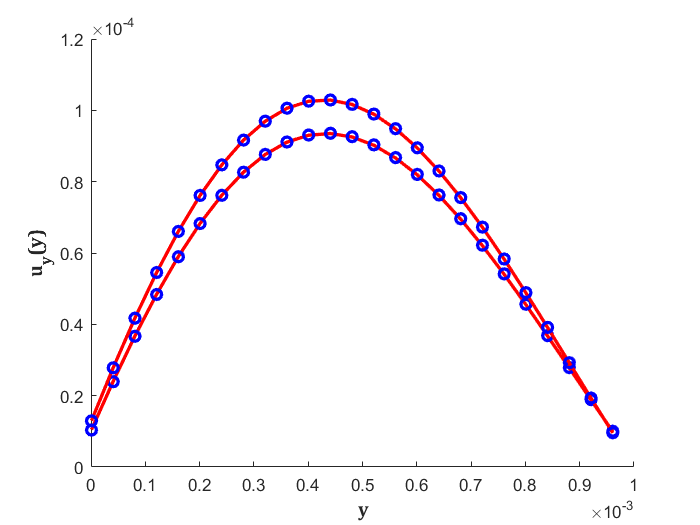}
\caption{Comparison of FEM solution $\mathbf{u}^h(t,x,y)$ with antisymmetric Lamb wave $\mathbf{u}(t,x,y)$, for $f_0= 600\,KHz$. Left: curves $u_x(\widetilde{t}_1,\overline{x},y)$ (solid red) and $u_x^h(\widetilde{t}_1,\overline{x},y)$ (blue o), right: curves $u_y(\widetilde{t}_1,\overline{x},y)$ (solid red) and $u_y^h(\widetilde{t}_1,\overline{x},y)$ (blue o). }
\label{isocurvas}
\end{figure}

In order to compare the FEM solution $\mathbf{u}^h(t,x,y)$ with $\mathbf{u}(t,x,y)$ for $f_0=600 \ KHz$, we fix $t=\widetilde{t}_1$, where  $\widetilde{t}_1=0.685 \cdot 10^{-5}$ is the arrival time at the point $\mathbf{p}_1$ on the top of the plate. The phase velocity $C(f_0)$ computed in the previous section and the zeros $t_1^{min},t_1^{max}$ of the function $u_y^h(t,\mathbf{p}_i)$ closest to $\widetilde{t}_1$ are also used to determine a suitable rectangle  $[x_1^{min},x_1^{max}] \times [0,L_y]$, where to compare $\mathbf{u}^h(\widetilde{t}_1,x,y)$ with $\mathbf{u}(\widetilde{t}_1,x,y)$. More precisely we compute $x_1^{min}=C(f_0)(\widetilde{t}_1-t_1^{min})=0.9 \cdot 10^{-2}$ and $x_1^{max}=C(f_0)(t_1^{max}-\widetilde{t}_1)=1.08 \cdot 10^{-2}$. In Figure \ref{isocurvas} we show the graphics of functions $\boldsymbol{u}^h(\widetilde{t}_1,\overline{x},y)$ and  $\boldsymbol{u}(\widetilde{t}_1,\overline{x},y)$ for $\overline{x}$ equal to $x_1^{min},x_1$ and $x_1^{max}$. As we observe, both components $u_x^h$ and $u_y^h$ of the approximated solutions curves $\boldsymbol{u}^h(\widetilde{t}_1,\overline{x},y)$ are very close to the corresponding components $u_x$ and $u_y$ of the antisymmetric Lamb wave $\boldsymbol{u}(\widetilde{t}_1,\overline{x},y)$.

\begin{figure}[h!]
\centering
\includegraphics[scale=0.22]{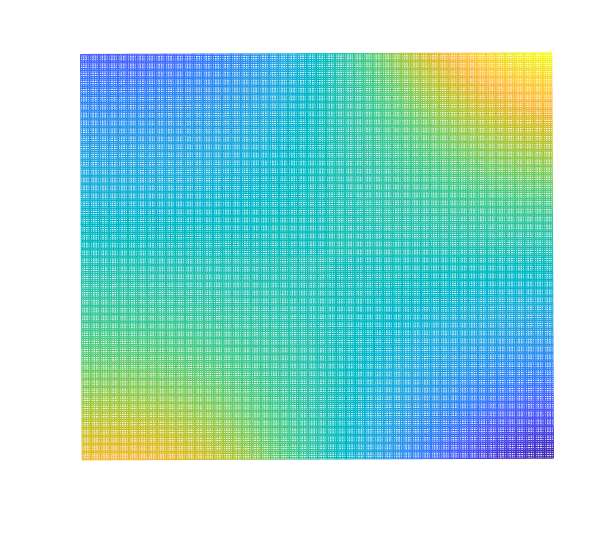}
\includegraphics[scale=0.235]{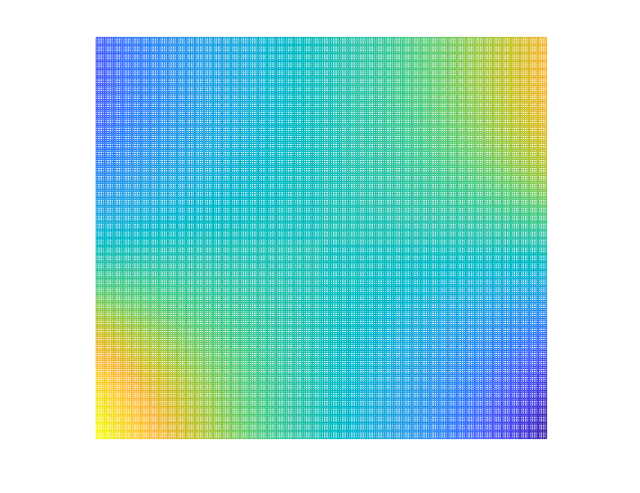}\\
\includegraphics[scale=0.21]{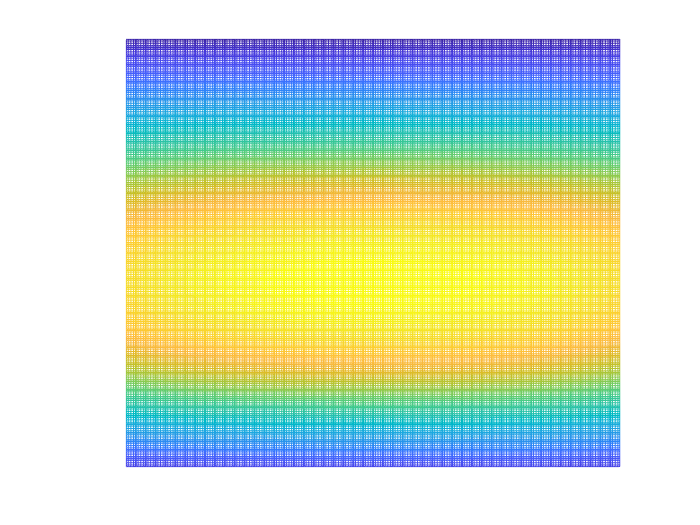}
\includegraphics[scale=0.21]{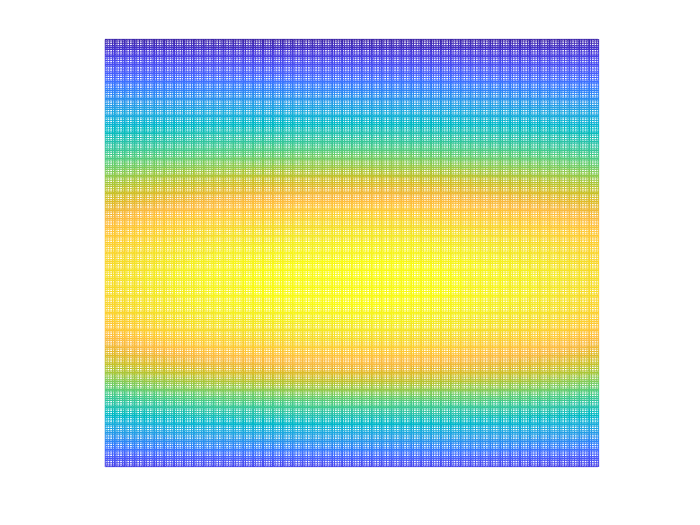}
\caption{First row: surfaces $u_x(\widetilde{t}_1,x,y)$ and $u_x^h(\widetilde{t}_1,x,y)$, second row: surfaces  $u_y(\widetilde{t}_1,x,y)$ and $u_y^h(\widetilde{t}_1,x,y)$, with $\widetilde{t}_1=0.685 \cdot 10^{-5}$, $(x,y) \in [0.9,1.08]\cdot 10^{-2} \times [0,1.0 \cdot 10^{-3}]$ and $f_0= 600\,KHz$.}
\label{superficies}
\end{figure}

To conclude this section, we show in Figure \ref{superficies} a 2D view of the surfaces $u^h_x(\widetilde{t}_1,x,y)$ and $u^h_y(\widetilde{t}_1,x,y)$ and compare them with the corresponding surfaces $u_x(\widetilde{t}_1,x,y)$ and $u_y(\widetilde{t}_1,x,y)$ for the antisymmetric Lamb wave. The comparison is  performed for $x \in [x_1^{min},x_1^{max}]=[0.9, 1.08]\cdot 10^{-2}$ and  $y \in [0,1.0 \cdot 10^{-3}]$. Here is also evident the good qualitative correspondence between the approximation $\boldsymbol{u}^h$ computed with \textit{quadratic} FEM and the Lamb wave solution $\boldsymbol{u}$.

\subsection{Study of convergence of  FEM approximation.}

In order to make a study of convergence, in this section we solve the propagation problem, for the pulse given by (\ref{g}) with parameters obtained for $f_0=600 \ KHz$ in Table \ref{tablaf}, and a sequence of  quasi-uniform and shape regular meshes $\{\tau_j\}_{j \geq 1}$ of size $h_j$, with  $h_j$ going to zero. The sequence $\{\tau_j\}_{j \geq 1}$ is obtained for an increasing sequence $\{ ny_j\}_{j \geq 1}$ of triangular vertices in the vertical direction. In terms of $ny_j$, the size  of the corresponding triangular mesh $ \tau_j$ is $h_j= \frac{\sqrt{2}}{ny_j }\cdot 10^{-3}\,$  and the number of degrees of freedom is equal to $100\,{{\it ny_j}}^{2}+102\,{\it ny_j}+2$, if Lagrange \textit{linear} finite elements are used, and equal to $400\,{{\it ny_j}}^{2}+204\,{\it ny_j}+2 $ in the case of Lagrange \textit{quadratic} finite elements.

The approximated solution $\boldsymbol{u}^h(t,x,y)$ is computed with FreeFem++. In Table \ref{tabla1} we show the behavior of the error $e_{\widetilde{t}}$,  between  FEM solutions corresponding to mesh sizes $h_{j-1}$ and $h_j$, for $\widetilde{t}=\widetilde{t}_1 = 0.685 \cdot 10^{-5}$. The error $e_{\widetilde{t}}$ is computed as,

$$e_{\widetilde{t}}\,\,:=\|\boldsymbol{u}^{h_{j-1}}(\widetilde{t},\boldsymbol{x})- \boldsymbol{u}^{h_j}(\widetilde{t},\boldsymbol{x})\| =\left(\int_{x^{min}}^{x^{max}}\int_{0}^{L_y} \left( \boldsymbol{u}^{h_{j-1}}(\widetilde{t},\boldsymbol{x})- \boldsymbol{u}^{h_{j}}(\widetilde{t},\boldsymbol{x})\right)^2 \;d\boldsymbol{x}\right)^{1/2} $$

where the values of $x^{min}=0.9 \cdot 10^{-2}$ and $x^{max}=1.08 \cdot 10^{-2}$ were computed using the procedure described in section \ref{secqual}. The corresponding number of degrees of freedom (dof), i.e. the dimension of the matrix $\boldsymbol{A}$, is also displayed for \textit{linear} and \textit{quadratic} FEM solutions.

\begin{table}[h!]
\caption{Error $e_{\widetilde{t}}$  between FEM solutions corresponding to mesh sizes $h_{j-1}$ and $h_j$,
for $\widetilde{t} = 0.685 \cdot 10^{-5}$. The solution $\boldsymbol{u}^h$
was computed using \textit{linear} and \textit{quadratic} Lagrange finite elements.}
\centering
\begin{tabular}{|c||c|c||c|c|}\hline
  & \multicolumn{2}{|c||}{Linear FEM } & \multicolumn{2}{|c|}{Quad. FEM } \\ \hline
$h_j$ & {\it dof} & $e_{\widetilde{t}}$ & {\it dof} & $e_{\widetilde{t}}$ \\ \hline
% $ h_j\,\,(m) $      &\# of freedom &$e_{\widetilde{t}}\,\,(m)$& \# of freedom  & $e_{\widetilde{t}}\,\,(m)$ \\ \hline
$ 7.07 \cdot 10^{-4}$  & $606 $      & $2.2768 \cdot 10^{-7} $     & $2010 $ & $7.8126 \cdot 10^{-8}  $    \\ \hline
$ 4.71 \cdot 10^{-4}$  & $1208$      & $1.2991 \cdot 10^{-7} $     & $4214 $ & $8.9191 \cdot 10^{-9}  $    \\ \hline
$ 3.53 \cdot 10^{-4}$  & $2010$      & $6.7139\cdot 10^{-8}  $     & $7218 $ & $2.0703 \cdot 10^{-9}  $    \\ \hline
$ 2.82 \cdot 10^{-4}$  & $3012$      & $3.6441\cdot 10^{-8}  $     & $11022$ & $7.3263 \cdot 10^{-10} $    \\ \hline
$ 2.35 \cdot 10^{-4}$  & $4214$      & $2.1537\cdot 10^{-8}  $     & $15626$ & $3.3669 \cdot 10^{-10} $    \\ \hline
$ 2.02 \cdot 10^{-4}$  & $5616$      & $1.3660\cdot 10^{-8}  $     & $21030$ & $1.8367 \cdot 10^{-10} $    \\ \hline
$ 1.76 \cdot 10^{-4}$  & $7218$      & $9.1647\cdot 10^{-9}  $     & $27234$ & $1.1266 \cdot 10^{-10} $    \\ \hline
$ 1.57 \cdot 10^{-4}$  & $9020$      & $6.4316\cdot 10^{-9}  $     & $34238$ & $7.5068 \cdot 10^{-11}  $    \\ \hline
$ 1.41 \cdot 10^{-4}$  & $11022$     & $4.6813\cdot 10^{-9}  $     & $42042$ & $5.3120 \cdot 10^{-11} $    \\ \hline
$ 1.28 \cdot 10^{-4}$  & $13224$     & $3.5115\cdot 10^{-9}  $     & $50646$ & $3.9344 \cdot 10^{-11} $    \\ \hline
$ 1.17 \cdot 10^{-4}$  & $15626$     & $2.7008\cdot 10^{-9}  $     & $60050$ & $3.0197 \cdot 10^{-11} $    \\ \hline
$ 1.08 \cdot 10^{-4}$  & $18228$     & $2.1213\cdot 10^{-9}  $     & $70254$ & $2.3837 \cdot 10^{-11} $    \\ \hline
$ 1.01 \cdot 10^{-4}$  & $21030$     & $1.6966\cdot 10^{-9}  $     & $81258$ & $1.9257 \cdot 10^{-11} $    \\ \hline
\end{tabular}
\label{tabla1}
\end{table}

\begin{figure}[h!]
\centering
\includegraphics[scale=0.4]{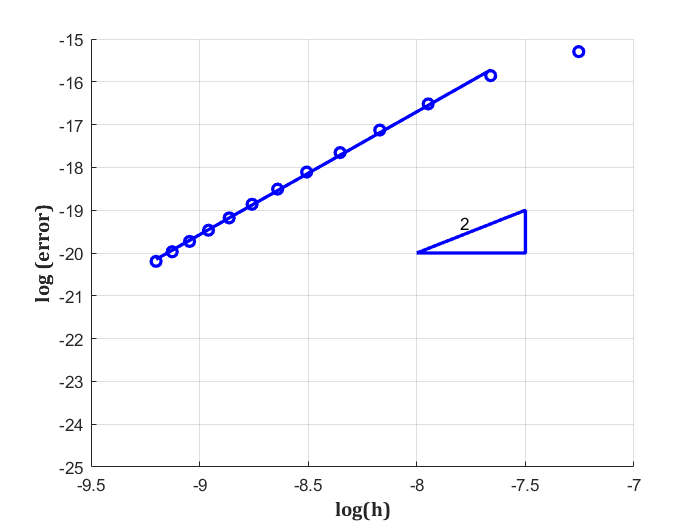}
\includegraphics[scale=0.4]{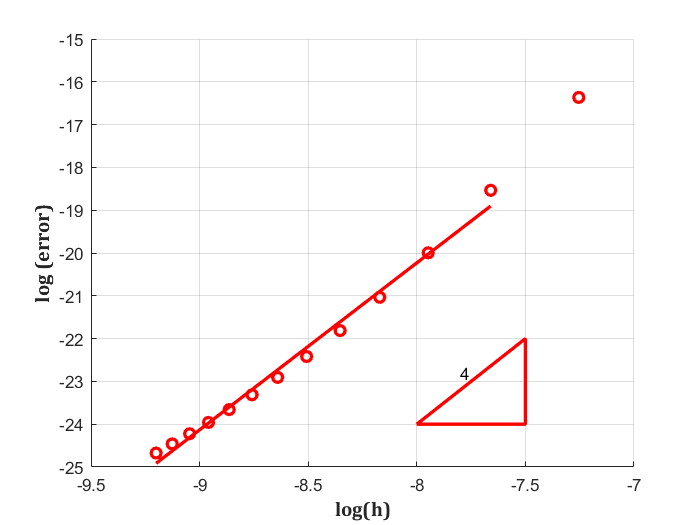}
\caption{Study of convergence of  FEM solutions for a sequence of meshes of size $h_j$, with  $h_j$ going to zero. Circles correspond to errors in Table \ref{tabla1}, continuous line is the best fitting line to errors, the graphic is in log-log scale.  Left: \textit{linear} FEM solution, right: \textit{quadratic} FEM solution.}
\label{convergencia}
\end{figure}

As we observe in Table \ref{tabla1} and Figure \ref{convergencia}, the error fulfills the \textit{a priori} error estimate in the $L_2$ norm proved in Theorem \ref{Theo3}, i.e., the error goes to zero with order $h^{k+1}$, where $k=1$  for Lagrange \textit{linear}  finite elements and $k=2$  for Lagrange \textit{quadratic} finite elements.

\section*{Conclusions}

Using finite differences  to approximate the temporal variable, the partial differential equations describing the propagation of an ultrasonic pulse along a thin plate is  approximately solved. For each fixed time, the corresponding problem in spacial variables is solved with classic Finite Element Method. Starting from the variational formulation of the problem, it is proved that the hypothesis of Lax-Milgram theorem holds and therefore the weak problem has an unique solution on a Hilbert space of functions. The Galerkin approach is used to compute the approximated solution on a finite dimensional space of piecewise polynomial functions defined on a triangulation of the physical domain. \textit{A priori} error estimates are obtained for a solution based on piecewise polynomials of degree $k$, showing that the approximate solution converges to the solution of the variational problem, when the size $h$ of the mesh goes to zero. Furthermore, it is proved that the energy norm of error is proportional to $h^{k}$, while under an additional stability hypothesis the $L_2$ norm of the error is proportional to $h^{k+1}$.

Numerical results were obtained with the software FreeFem++, using Lagrange \textit{linear} and \textit{quadratic} finite elements. The approximated solution is compared with the analytical solution of a similar wave propagation problem for fixed values of the temporal variable, and the good qualitative correspondence between the displacements in both wave propagation problems becomes apparent. In a convergence study, the approximated solution is computed for decreasing mesh sizes and it is shown  that the error decreases with the theoretically expected velocity as the mesh size tends to zero.

Moreover, we develop a successful strategy to compute points on the phase velocity dispersion curve. This strategy could be also used in another problems, where no analytical solution exists, such as in the case of more complicated geometries or  nonelastic and anisotropic materials. Due the importance of the dispersion curves for industrial applications, this subject will be treated in a future research.

\bibliographystyle{plain}
\bibliography{ManuelBibl}

\end{document}